\documentclass{amsart}

\usepackage{amsmath}
\usepackage{amsthm}
\usepackage{amssymb}
\usepackage{ifthen}

\input diagrams

\theoremstyle{plain}
\newtheorem*{main}{Main~Theorem}
\newtheorem{theo}{Theorem}[section]
\newtheorem{lemm}[theo]{Lemma}
\newtheorem{prop}[theo]{Proposition}
\newtheorem{coro}[theo]{Corollary}

\theoremstyle{definition}

\newtheorem{defi}[theo]{Definition}
\newtheorem{rema}[theo]{Remark}

\newcommand{\bq}{\begin{equation}}

\newcommand{\enbe}{\end{equation}}

\renewcommand{\a}{\alpha}
\renewcommand{\b}{\beta}

\newcommand{\gA}{{\mathfrak{A}}}
\newcommand{\ZZZ}{{\mathbf{Z}}}
\newcommand{\zz}{{\mathbb{Z}}}

\renewcommand{\leq}{\leqslant}
\renewcommand{\geq}{\geqslant}

\newcommand{\Id}{\text{\rm Id}}
\newcommand{\id}{\text{id}}
\newcommand{\qs}{\hfill\square}
\newcommand{\enre}{\end{rema}}
\newcommand{\enco}{\end{coro}}
\newcommand{\enpr}{\end{prop}}
\newcommand{\enth}{\end{theo}}
\newcommand{\enle}{\end{lemm}}
\newcommand{\enen}{\end{enumerate}}
\newcommand{\eneq}{\end{equation}}
\newcommand{\enal}{\end{aligned}}

\renewcommand{\beth}{\begin{theo}}
\newcommand{\bele}{\begin{lemm}}
\newcommand{\bepr}{\begin{prop}}
\newcommand{\beeq}{\begin{equation}}
\newcommand{\bega}{\begin{gather}}
\newcommand{\been}{\begin{enumerate}}

\newcommand{\arar}{A_\rho}
\newcommand{\arz}{\arar [[z]]}

\newcommand{\arzz}{\arar ((z))}

\newcommand{\woar}{W_1(A,\rho)}
\newcommand{\war}{W(A,\rho)}
\newcommand{\Nil}{{\rm Nil}}
\newcommand{\warab}{W(A,\rho)^{ab}}

\newcommand{\karzz}{K_1(\arzz)}

\newcommand{\nil}{{\rm Nil}_0(A,\rho)}

\newcommand{\bhs}
{H.Bass, A.Heller, R.G.Swan,
\emph{The Whitehead group of a polynomial extension},
Inst. Hautes Etudes Sci. Publ. Math. {\bf 22} (1964),
61--79
}

\newcommand{\farhsi}
{F.T.Farrell, W.-C.Hsiang,
\emph{A formula for $K_1R_\alpha[T]$},
Proc. Symp. Pure Math., Vol. {\bf 17} (1968), 192--218}

\newcommand{\patou}
{ A.V.Pajitnov, \emph{ On the Novikov complex for rational Morse forms},
Annales de la Facult\'e de Sciences de Toulouse {\bf 4}  (1995), 297--338.
}

\newcommand{\bass}
{ H.Bass,
\emph{Algebraic K-theory},
Benjamin, 1968.
}

\newcommand{\panonab}{ A.V.Pajitnov,
\emph{Closed orbits of gradient flows
and logarithms of non-abelian Witt vectors
},
e-print math.DG/9908010, to appear in $K$-theory}

\newcommand{\ranibook}
{   A.A.Ranicki,
\emph{High-dimensional knot theory,         }
Springer, 1998
}

\newcommand{\ranitor}
{   A.A.Ranicki,
\emph{The algebraic theory of torsion I.,}
Lecture Notes in Mathematics {\bf 1126} (1985), 199--237,
Springer.
}

\newcommand{\sieben}
{L.Siebenmann,
\emph{A total Whitehead torsion obstruction to fibering over the circle},
Comment. Math. Helv. {\bf 45} (1970), 1--48.}

\begin{document}

\title[The Whitehead group of the Novikov ring]
{The Whitehead group of the Novikov ring}
\author{A.V.Pajitnov}
\address[A.V.Pajitnov]{Universit\'e de Nantes\newline
\indent UMR 6629 CNRS, D\'epartement de Math\'ematiques\newline
\indent 2, rue de la Houssini\`ere\newline
\indent 44072 Nantes Cedex\newline
\indent France}
\email{pajitnov@math.univ-nantes.fr}
\author{A.A.Ranicki}
\address[A.A.Ranicki]{Dept. of Mathematics and Statistics\newline
\indent University of Edinburgh\newline
\indent Edinburgh EH9 3JZ\newline
\indent Scotland, UK}
\email{aar@maths.ed.ac.uk}
\maketitle

\section*{Abstract}

The Bass-Heller-Swan-Farrell-Hsiang-Siebenmann decomposition of the
Whitehead group $K_1(A_{\rho}[z,z^{-1}])$
of a twisted Laurent polynomial extension $A_{\rho}[z,z^{-1}]$
of a ring $A$ is generalized to a decomposition of the Whitehead
group $K_1(A_{\rho}((z)))$ of a twisted Novikov ring of power series
$A_{\rho}((z))=A_{\rho}[[z]][z^{-1}]$. The decomposition
involves a summand $W_1(A,\rho)$ which is an abelian
quotient of the multiplicative group $W(A,\rho)$ of Witt vectors
$1+a_1z+a_2z^2+\dots \in A_{\rho}[[z]]$. An example is constructed
to show that in general the natural surjection $W(A,\rho)^{ab} \to W_1(A,\rho)$
is not an isomorphism.

\section*{Introduction }

We obtain a splitting theorem for the Whitehead group $K_1(A_{\rho}((z)))$
of the Novikov ring of power series $A_{\rho}((z))$, which
is an analogue of the well-known splitting theorem for the Whitehead
group $K_1(A_{\rho}[z,z^{-1}])$ of the Laurent ring of polynomials $A_{\rho}[z,z^{-1}]$.

Let $A$ be an associative ring with 1.  Given an automorphism $\rho:A
\to A$ let $z$ be an indeterminate over $A$ such that
$$az~=~z\rho(a)~~(a \in A)~.$$
The {\it $\rho$-twisted polynomial extension ring} $A_{\rho}[z]$ is the
ring of polynomials $\sum\limits^{\infty}_{j=0}a_jz^j$
with $a_j=0\in A$ for all but a finite number of $j\geq 0$.
The {\it $\rho$-twisted power series ring} $A_{\rho}[[z]]$ is the
ring of power series $\sum\limits^{\infty}_{j=0}a_jz^j$ for arbitrary
$a_j \in A$.
Inverting $z$ we also obtain the {\it $\rho$-twisted Laurent
polynomial extension ring} $A_{\rho}[z,z^{-1}]$
of polynomials $\sum\limits^{\infty}_{j=-\infty}a_jz^j$
with $a_j=0\in A$ for all but a finite number of $j\in \ZZZ$, and
the {\it $\rho$-twisted Novikov formal power series ring} $A_{\rho}((z))$
of polynomials $\sum\limits^{\infty}_{j=-\infty}a_jz^j$
with $a_j=0\in A$ for all but a finite number of $j<0$.

The Whitehead groups of the polynomial rings split as
$$\begin{array}{l}
K_1(A_{\rho}[z])~=~K_1(A)\oplus \widetilde{\text{Nil}}_0(A,\rho)~,\\[1ex]
K_1(A_{\rho}[z,z^{-1}])~=~K_1(A,\rho)\oplus
\widetilde{\text{Nil}}_0(A,\rho)\oplus \widetilde{\text{Nil}}_0(A,\rho^{-1})
\end{array}$$
with $K_1(A,\rho)$ the class group of pairs $(P,\phi)$ with $P$ a
f.g. (= finitely generated) projective $A$-module and $\phi:P \to P$ a
$\rho$-twisted
automorphism, and $\widetilde{\text{Nil}}_0(A,\rho)$ the reduced class
group of pairs $(P,\nu)$ with $P$ a f.g. projective $A$-module and
$\nu:P \to P$ a nilpotent $\rho$-twisted endomorphism (Bass, Heller and
Swan \cite{bhs}, Bass \cite{bass}, Farrell and Hsiang \cite{farhsi},
Siebenmann \cite{sieben}).

The augmentation
$$A_{\rho}[[z]] \to A~;~\sum\limits^{\infty}_{j=0}a_jz^j \mapsto a_0$$
induces a split surjection $K_1(A_{\rho}[[z]]) \to K_1(A)$, with the kernel
$$NK_1(A_{\rho}[[z]])~=~{\rm ker}(K_1(A_{\rho}[[z]]) \to K_1(A))$$
such that
$$K_1(A_{\rho}[[z]])~=~K_1(A) \oplus NK_1(A_{\rho}[[z]])~.$$
Pajitnov \cite{patou} identified $NK_1(A_{\rho}[[z]])$ with the
subgroup $W_1(A,\rho) \subseteq K_1(A_{\rho}[[z]])$ represented
by the Witt vectors, that is the units in $A_{\rho}[[z]]$ of the type
$$w~=~1 +\sum\limits^{\infty}_{j=1}a_jz^j \in A_{\rho}[[z]]^{\bullet}~.$$

\begin{main} \label{t:main} The Whitehead group of the Novikov ring splits as
$$K_1(A_{\rho}((z)))~=~K_1(A,\rho)\oplus W_1(A,\rho) \oplus
\widetilde{\mbox{\rm Nil}}_0(A,\rho^{-1})~.$$
\end{main}

This splitting was obtained in the untwisted case $\rho=1$ in Chapter
14 of Ranicki \cite{ranibook} by using general results on the algebraic
$K$-theory of localization-completion squares, such as
\begin{diagram}
A_{\rho}[z] & \rTo  & A_{\rho}[z,z^{-1}] \\
   \dTo     &       &  \dTo   \\
A_{\rho}[[z]] &  \rTo   & A_{\rho}((z))  \\
\end{diagram}

\noindent
In principle, the general method also works in the twisted case, using the
equivalences of exact categories
$$\matrix
\{\hbox{\rm f.g. projective $A$-modules $P$ with a nilpotent $\rho$-twisted
endomorphism $\nu:P \to P$}\}\vspace{1.5mm}\\
\approx
\{\hbox{\rm $z$-primary torsion $A_{\rho}[z]$-modules of homological
dimension 1}\}\vspace{1.5mm}\\
\hspace{2mm}
\approx \{\hbox{\rm $z$-primary torsion $A_{\rho}[[z]]$-modules of
homological dimension 1}\}
\endmatrix$$
but here we prefer to use a direct method.
Chapter \ref{s:main} of this paper
is devoted to a direct proof of the Main Theorem, which follows the
direct proof of the splitting theorem for $K_1(A_\rho[z,z^{-1}])$,
except for the (Higman) linearization result which does not have an analogue
for $K_1(A_{\rho}((z)))$.

Our Main Theorem is used in the work of the first author \cite{panonab}
to define a non-abelian logarithm in the twisted case.  Here is the
corollary which is used in \cite{panonab}, and which follows
immediately from the proof of the Main Theorem, given in Section
\ref{su:whnovikov}.

\begin{coro}\label{c:for_zeta}
The homomorphism $\widehat C_2:W_1(A,\rho)\to\karzz$ induced by the inclusion
has a left inverse $\widehat B_2:\karzz\to W_1(A,\rho)$, which satisfies:
\been
\item $\widehat B_2$ vanishes on the image of $K_1(A)$ in $\karzz$,
\item $\widehat B_2(\tau(z))=0$, with $\tau(z) \in \karzz$ the
torsion of the invertible $1\times  1$-matrix $(z)$.\hfill $\qs$
\enen
\end{coro}

Write the multiplicative group of Witt vectors in $A_{\rho}[[z]]$ as
$$W(A,\rho)~=~1+zA_{\rho}[[z]] \subseteq A_{\rho}[[z]]^{\bullet}~,$$
with abelianization
$$W(A,\rho)^{ab}~=~W(A,\rho)/[W(A,\rho),W(A,\rho)]~.$$
It was claimed in Proposition 14.6 of Ranicki \cite{ranibook} that the
natural surjection $W(A,\rho)^{ab} \to W_1(A,\rho)$ is an injection, at
least in the untwisted case $\rho=1$.  In Chapter \ref{s:ex} of this paper we
correct this, constructing a family of explicit counterexamples,
already in the untwisted case $\rho=1$.

\subsection{Acknowledgements}
\label{su:ack}

The main part of this work was done in November 1999.  We had an
opportunity to discuss the Whitehead groups of Novikov rings when the
first author stayed in Great Britain by the invitation based on the
EPRSC grant GR/M98159.  A.Pajitnov would like to express
here his gratitude to the Mathematical Institute of the Oxford
University for hospitality and financial support during this visit.

\section{Class and torsion}
\label{s:class}

We recall the definitions of the $K$-groups $K_0,K_1$ of additive
and exact categories, and also of the less familiar isomorphism torsion
group of Ranicki \cite{ranitor}.

The {\it class group} $K_0({\mathbb C})$ of an exact category ${\mathbb C}$
is the abelian group with one generator $[M]$ for
each isomorphism class of objects $M$ in ${\mathbb C}$, and one relation
$$[L]-[M]+[N]~=~0$$
for each short exact sequence in ${\mathbb C}$
$$0 \to L \to M \to N \to 0~.$$
\indent
The {\it Whitehead group} $K_1({\mathbb C})$ of an exact category
${\mathbb C}$ has one generator $\tau(f)$ for each automorphism
$f:M \to M$ in ${\mathbb C}$, and relations
$$\aligned
{}&{\rm (i)}~\tau(e)-\tau(f)+\tau(g)~=~0\\
{}&\hskip10mm
\hbox{\rm for each automorphism of a short exact sequence in ${\mathbb C}$}\\[2ex]
{}&\hskip10mm
\begin{diagram}
0 & \rTo & L & \rTo  & M & \rTo  & N & \rTo  & 0 \\
       &  & \dTo^e   &       &  \dTo^f   &       &  \dTo^g  & \\
0 & \rTo & L & \rTo  & M & \rTo  & N & \rTo  & 0
\end{diagram}\\[2ex]
&{\rm (ii)}~\tau(gf:M \to M)~=~\tau(f:M\to M) + \tau(g:M \to M)\\
&\hskip10mm \hbox{\rm for automorphisms $f,g:M \to M$ in ${\mathbb C}$}~.
\endaligned$$

The algebraic $K$-groups $K_0(A)$, $K_1(A)$ of a ring $A$
are defined by
$$K_i(A)~=~K_i({\mathbb P}(A))$$
with ${\mathbb P}(A)$ the exact category of f.g. projective $A$-modules;
$K_1(A)$ is called the Whitehead group of $A$.

\medskip
\bepr\label{Proposition1}
{\rm (Bass \cite{bass},p.\ 397, Ranicki \cite{ranitor},Prop.\ 1.1)}\\
Let $\mathbb A$ be an additive category with the split exact structure\,:
a sequence
\begin{diagram}
0 & \rTo & L & \rTo^i & M & \rTo^j  & N & \rTo  & 0
\end{diagram}
is exact if and only if there exists a morphism $u:N \to M$ such that
$$(i~u)~:~L \oplus N \to M$$
is an isomorphism.\\
{\rm (i)} The Whitehead group $K_1({\mathbb A})$ is (isomorphic to) the abelian
group $K'_1({\mathbb A})$ with one generator $\tau(f)$
for each automorphism $f:M \to M$ in ${\mathbb A}$, and relations
$$\aligned
{}&\tau(\hbox{$\begin{pmatrix} f & d \\ 0 & f' \end{pmatrix}$}:
M \oplus M'\to M\oplus M')~=~\tau(f:M \to M) + \tau(f':M'\to M')\\
{}&\hbox{\it
for automorphisms $f:M \to M$, $f':M'\to M'$ and any morphism $d:M' \to M$},\\
{}&\tau(gf:M \to M)~=~\tau(fg:N \to N)~
\hbox{\it for isomorphisms $f:M \to N$, $g:N \to M$ in ${\mathbb A}$}~.
\endaligned$$
{\rm (ii)} The Whitehead group $K_1({\mathbb A})$ is (isomorphic to) the abelian
group $K''_1({\mathbb A})$ with one generator $\tau(f)$
for each automorphism $f:M \to M$ in ${\mathbb A}$, and relations
$$\aligned
{}&\tau(f\oplus f':M\oplus M' \to M\oplus M')~=~
\tau(f:M \to M)+\tau(f':M' \to M')\\
{}&\hskip30mm
(\hbox{\it for automorphisms $f:M \to M$, $f':M' \to M'$})~,\\
{}&\tau(gf:M \to M)~=~\tau(f:M \to M)+\tau(g:M \to M)\\
{}&\hskip30mm
(\hbox{\it for automorphisms $f:M \to M$, $g:M \to M$})~,\\
{}&\tau(gf:M \to M)~=~\tau(fg:N \to N)\\
{}&\hskip30mm (\hbox{\it for isomorphisms $f:M \to N$, $g:N \to M$})~.
\endaligned$$
\enpr
\noindent{\it Proof.}
Let $S$ be the set of automorphisms in ${\mathbb A}$, and let
$R,R',R'' \subset {\mathbb Z}[S]$ be the subgroups defined by
$$\begin{array}{l}
R~=~\{~\tau(e)-\tau(f)+\tau(g)~\hbox{\rm for an exact sequence}~
0\to e\to f\to g\to 0~,\\[1ex]
\hphantom{R~=~\{~} \tau(gf)-\tau(f)-\tau(g)~{\rm for~automorphisms}~f,g\}~,\\[1ex]
R'~=~\{~\tau(f)-
\tau\hbox{$\begin{pmatrix} f & d \\ 0 & f' \end{pmatrix}$}+\tau(f')~
{\rm for~automorphisms}~f,f',\\
\hphantom{R'~=~\{~}
\tau(gf)-\tau(fg)~{\rm for~isomorphisms}~f,g~\}~,\\[1ex]
R''~=~\{~\tau(f)-\tau(f\oplus f')+\tau(f')~{\rm for~automorphisms}~f,f'~,\\[1ex]
\hphantom{R''~=~\{~}
\tau(f)-\tau(gf)+\tau(g)~{\rm for~automorphisms}~f,g~,\\[1ex]
\hphantom{R''~=~\{~}
\tau(gf)-\tau(fg)~{\rm for~isomorphisms}~f,g~\}~,
\end{array}$$
so that
$$K_1({\mathbb A})~=~{\mathbb Z}[S]/R~~,~~
K'_1({\mathbb A})~=~{\mathbb Z}[S]/R'~~,~~
K''_1({\mathbb A})~=~{\mathbb Z}[S]/R''~.$$
We shall prove that $R=R'=R''$ by first showing that
$R \subseteq R' \subseteq R$ and then
$R' \subseteq R''\subseteq R$.\newline
(i) ($R \subseteq R'$)
Given an automorphism of a short exact sequence
in ${\mathbb A}$
\begin{diagram}
0 & \rTo & L & \rTo^i  & M & \rTo^j  & N & \rTo  & 0 \\
  &    &  \dTo^e   &       &  \dTo^f   &       &  \dTo^g  & \\
0 & \rTo & L & \rTo^i  & M & \rTo^j  & N & \rTo  & 0
\end{diagram}
there exists a morphism $u:N \to M$ such that
$$h~=~(i~u)~:~L \oplus N \to M$$
is an isomorphism, with inverse of the form
$$h^{-1}~=~\begin{pmatrix} v \\ j \end{pmatrix}~:~M \to L \oplus N~.$$
Now
$$h^{-1}fh~=~
\begin{pmatrix} e  &  vfu \\ 0 & g \end{pmatrix}~:~L\oplus N \to L \oplus N~.$$
Thus
$$\tau(e)-\tau(f)+\tau(g)~=~(\tau(e)-\tau(h^{-1}fh)+\tau(g))+
(\tau(h^{-1}fh)-\tau(f)) \in R'$$
and $R \subseteq R'$.\newline
(ii) ($R' \subseteq R$) Every automorphism in $\mathbb A$ is of the form
$$\alpha~=~\begin{pmatrix} f & d \\ 0 & f' \end{pmatrix}~:~
M \oplus M'\to M\oplus M'$$
and fits into an automorphism of a short exact sequence
\begin{diagram}
0 &\rTo & M & \rTo  & M\oplus M' & \rTo  & M' & \rTo  & 0 \\
  &     &  \dTo^f   &       &  \dTo^{\alpha}   &       &  \dTo^{f'}  & \\
0 &\rTo & M & \rTo  & M\oplus M'& \rTo  & M' & \rTo  & 0
\end{diagram}
so that
$$\tau(\alpha)-\tau(f)-\tau(f') \in R~.$$
(iii) ($R' \subseteq R''$)
For any elementary automorphism in $\mathbb A$
$$\hbox{$\begin{pmatrix} 1 & e \\ 0 & 1 \end{pmatrix}$}:
M \oplus M'\to M\oplus M'$$
the automorphisms
$$\begin{array}{l}
g~=~\hbox{$\begin{pmatrix} 1 & 0 & 1 \\ 0 & 1 & 0 \\ 0 & 0 & 1 \end{pmatrix}$}~:~
M \oplus M' \oplus M' \to M \oplus M' \oplus M'~,\\[4ex]
h~=~\hbox{$\begin{pmatrix} 1 & 0 & 0 \\ 0 & 1 & 0 \\ 0 & e & 1 \end{pmatrix}$}~:~
M \oplus M' \oplus M' \to M \oplus M' \oplus M'
\end{array}$$
are such that
$$\hbox{$\begin{pmatrix} 1 & e \\ 0 & 1 \end{pmatrix}$}\oplus  1~=~
ghg^{-1}h^{-1}~:~M \oplus M' \oplus M' \to M \oplus M' \oplus M'$$
so that
$$\tau\hbox{$\begin{pmatrix} 1 & e \\ 0 & 1 \end{pmatrix}$} \in R''~.$$
More generally, for any automorphism of the type
$$\hbox{$\begin{pmatrix} f & e \\ 0 & f' \end{pmatrix}$}:
M \oplus M'\to M\oplus M'$$
we have
$$\hbox{$\begin{pmatrix} f & 0 \\ 0 & f' \end{pmatrix}$}^{-1}
\hbox{$\begin{pmatrix} f & e \\ 0 & f' \end{pmatrix}$}~=~
\hbox{$\begin{pmatrix} 1 & f^{-1}e \\ 0 & 1 \end{pmatrix}$}$$
so that
$$\tau\hbox{$\begin{pmatrix} f & e \\ 0 & f' \end{pmatrix}$}-
\tau\hbox{$\begin{pmatrix} f & 0 \\ 0 & f' \end{pmatrix}$} \in R''~.$$
(iv) ($R'' \subseteq R'$)
For any automorphisms $f:M \to M$, $g:M \to M$
$$\begin{pmatrix} gf & 0 \\ 0 & 1 \end{pmatrix}~=~
\begin{pmatrix} g & 0\\ 0 & f \end{pmatrix}
\begin{pmatrix} f & 0 \\ 0 & f^{-1} \end{pmatrix}~:~
M \oplus M \to M \oplus M$$
and for any automorphism $\alpha:M \to M$ we have the Whitehead lemma
identity
$$\begin{pmatrix} \alpha & 0 \\ 0 & \alpha^{-1} \end{pmatrix}~=~
\begin{pmatrix} 1 & \alpha \\ 0 & 1 \end{pmatrix}
\begin{pmatrix} 1 & 0 \\ -\alpha^{-1} & 1 \end{pmatrix}
\begin{pmatrix} 1 & \alpha \\ 0 & 1 \end{pmatrix}
\begin{pmatrix} 0 & -1 \\ 1 & 0 \end{pmatrix}~:~M \oplus M \to M \oplus M~,$$
so that
$$\tau(f) - \tau(gf) +\tau(g) \in R'~.$$
For any isomorphisms $f:M \to N$, $g:N \to M$
$$\begin{pmatrix} gf & 0 \\ 0 & (fg)^{-1} \end{pmatrix}~=~
\begin{pmatrix} 0 & g\\ -g^{-1} & 0 \end{pmatrix}
\begin{pmatrix} 0 & -f^{-1} \\ f & 0 \end{pmatrix}~:~M \oplus N \to M \oplus N~,$$
and for any isomorphism $h:M \to N$ we have
$$\begin{pmatrix} 0 & -h^{-1} \\ h & 0 \end{pmatrix}~=~
\begin{pmatrix} 1 & 0 \\ h & 1 \end{pmatrix}
\begin{pmatrix} 1 & -h^{-1} \\ 0 & 1 \end{pmatrix}
\begin{pmatrix} 1 & 0 \\ h & 1 \end{pmatrix}~:~
M \oplus N \to M \oplus N$$
so that
$$\tau(gf) - \tau(fg) \in R'~.$$
Thus $R'' \subseteq R'$.\hfill$\qs$

\begin{defi}\label{Definition1} (Ranicki \cite{ranitor})
The {\it isomorphism torsion group} $K^{iso}_1({\mathbb A})$
of an additive category $\mathbb A$ is the abelian group with one generator
$\tau^{iso}(f)$ for each isomorphism $f:M \to N$ in $\mathbb A$, and
relations
$$\aligned
{}&\tau^{iso}(gf:M \to P)~=~\tau^{iso}(f:M \to N)+\tau^{iso}(g:N \to P)~,\\
{}&\tau^{iso}(f\oplus f':M\oplus M' \to N\oplus N')~=~
\tau^{iso}(f:M \to N)+\tau^{iso}(f':M' \to N')~.
\endaligned$$\hfill$\qs$
\end{defi}

Every automorphism is an isomorphism, so there is an evident forgetful map
$$K_1({\mathbb A}) \to K^{iso}_1({\mathbb A})~;~\tau(f) \mapsto \tau^{iso}(f)~.$$
This map is an isomorphism if every isomorphism in ${\mathbb A}$ is an
automorphism.

A functor $F:{\mathbb A} \to {\mathbb B}$ of additive categories induces
morphisms of the torsion groups
$$\aligned
{}&F~:~K_1({\mathbb A}) \to K_1({\mathbb B})~;~
\tau(f:M \to M) \mapsto \tau(F(f):F(M) \to F(M))~,\\
{}&F~:~K^{iso}_1({\mathbb A}) \to K^{iso}_1({\mathbb B})~;~
\tau^{iso}(f:M \to N) \mapsto \tau^{iso}(F(f):F(M) \to F(N))~.
\endaligned$$

\begin{rema}\label{Remark1}
(i) There is an essential difference between
$K_1({\mathbb A})$ and $K^{iso}_1({\mathbb A})$. An equivalence of additive
categories $F:{\mathbb A} \to {\mathbb B}$ induces
an isomorphism $F:K_1({\mathbb A}) \cong K_1({\mathbb B})$, but the
induced morphism $F:K^{iso}_1({\mathbb A}) \to K^{iso}_1({\mathbb B})$
may not be an isomorphism. See (ii) below for an example.\\
(ii)  For any ring $A$ let
${\mathbb B}={\mathbb B}(A)$ be the additive category with objects based
f.g. free $A$-modules and $A$-module morphisms, and let ${\mathbb A}\subset
{\mathbb B}$ be the full subcategory with objects $A^n$ with the standard
basis $e=\{e_1,e_2,\dots,e_n\}$, where  $e_i=(0,\dots,0,1,0,\dots,0)$.
The inclusion $F:{\mathbb A} \to {\mathbb B}$ is an equivalence of
categories such that
$F:K^{iso}_1({\mathbb A})=K_1(A) \to K^{iso}_1({\mathbb B})$
is not an isomorphism: if $b=\{b_1,b_2,\dots,b_n\}$ is a non-standard
basis for $A^n$ and $f:(A^n,e) \to (A^n,b)$ is the isomorphism in ${\mathbb B}$
defined by $f=1:A^n \to A^n$ then
$\tau^{iso}(f) \in K^{iso}_1({\mathbb B})$ is not in the image of $F$.
\hfill$\qs$
\end{rema}

\medskip

\begin{defi}\label{Definition3}
(i) Given a category ${\mathbb A}$ let ${\rm Iso}({\mathbb A})$
be the subcategory with the same objects, but only the isomorphisms as morphisms.\\
(ii) An {\it isomorphism torsion structure} $F$ on an additive category
$\mathbb A$ is a functor
$$F~:~{\rm Iso}({\mathbb A}) \to {\rm Iso}({\mathbb A})$$
which is the identity on objects, and such that
$$\begin{array}{l}
F(f \oplus f')~=~F(f) \oplus F(f')~:\\[1ex]
\hphantom{F(f \oplus f')~=~}
F(M \oplus M')~=~M \oplus M' \to  F(N \oplus N')~=~N \oplus N'
\end{array}$$
for any isomorphisms $f:M \to N$, $f':M' \to N'$.\\
(iii) The {\it $F$-relative torsion} of an isomorphism $f:M \to N$ in $\mathbb A$ is
$$\tau^F(f)~=~\tau(F(f)^{-1}f:M \to M) \in K_1({\mathbb A})~.$$
(iv) An isomorphism torsion structure $F$ is {\it idempotent} if $F^2=F$.
\hfill$\qs$
\end{defi}
\medskip

\bepr\label{Proposition2}
Let ${\mathbb A}$ be an additive category, and let
$F:{\rm Iso}({\mathbb A}) \to {\rm Iso}({\mathbb A})$ be an isomorphism torsion structure.\\
{\rm (i)} The $F$-relative torsion function defines a morphism
$$\tau^F~:~K^{iso}_1({\mathbb A}) \to K_1({\mathbb A})~;~
\tau^{iso}(f:M \to N) \mapsto \tau^F(f)~=~\tau(F(f)^{-1}f:M \to M)~.$$
{\rm (ii)} $F$ induces an endomorphism of the Whitehead group
$$F~:~K_1({\mathbb A}) \to K_1({\mathbb A})~;~\tau(f:M \to M) \to
\tau(F(f):M \to M)$$
such that
$$1-F~:~K_1({\mathbb A}) \to K^{iso}_1({\mathbb A})
\rTo^{\tau^F} K_1({\mathbb A})~.$$
{\rm (iii)} If $F$ is idempotent
$$\aligned
{}&F^2~=~F~:~K_1({\mathbb A}) \to K_1({\mathbb A})~,\\
{}&{\rm im}(\tau^F:K^{iso}_1({\mathbb A}) \to K_1({\mathbb A}))~
=~{\rm im}(1-F:K_1({\mathbb A}) \to K_1({\mathbb A}))\\
{}&\hphantom{{\rm im}(\tau^F:K^{iso}_1({\mathbb A}) \to K_1({\mathbb A}))~}
=~{\rm ker}(F:K_1({\mathbb A}) \to K_1({\mathbb A}))~,\\
{}&K_1({\mathbb A}) ~=~
{\rm im}(F:K_1({\mathbb A}) \to K_1({\mathbb A}))\oplus
{\rm im}(\tau^F)~.
\endaligned$$
\enpr
\noindent{\it Proof.} Immediate from \ref{Proposition1}.\hfill$\qs$
\medskip

\begin{defi} \label{Definition4}
Let $A$ be a ring such that the rank of f.g. free $A$-modules is well-defined.
The {\it canonical idempotent isomorphism torsion structure} on the
additive category ${\mathbb B}(A)$ of based f.g. free $A$-modules
$$F^{can}~:~{\rm Iso}({\mathbb B}(A)) \to {\rm Iso}({\mathbb B}(A))$$
sends every isomorphism $f:M \to N$ to the isomorphism $F^{can}(f):M \to N$
sending the given basis of $M$ to the given basis of $N$.\hfill$\qs$
\end{defi}
\medskip

\bepr\label{Proposition3}
{\rm (i)} The canonical idempotent isomorphism torsion structure $F^{can}$
on ${\mathbb B}(A)$ induces
$$F^{can}~=~0~:~K_1({\mathbb B}(A))~=~K_1(A) \to K_1(A)~.$$
{\rm (ii)} The $F^{can}$-relative torsion
$$\begin{array}{l}
\tau^{F^{can}}~:~K^{iso}_1({\mathbb B}(A)) \to K_1({\mathbb B}(A))~=~K_1(A)~;\\[1ex]
\hphantom{\tau^{F^{can}}~:~}
\tau^{iso}(f:M \to N) \mapsto \tau^{F^{can}}(f)~=~
\tau(F^{can}(f)^{-1}f:M \to M)
\end{array}$$
is the standard way of assigning a torsion to an isomorphism of based
f.g. free $A$-modules, defining a surjection splitting the forgetful map
$K_1({\mathbb B}(A)) \to K^{iso}_1({\mathbb B}(A))$.
{\rm (iii)} The $F^{can}$-relative torsion of a chain equivalence $f:D \to E$ of finite
chain complexes in ${\mathbb B}(A)$ is the usual torsion $\tau(f) \in K_1(A)$.\hfill$\qs$
\enpr
\medskip

Given a ring morphism $F:A \to B$ regard $B$ as a $(B,A)$-bimodule by
$$B \times B \times A \to B~;~(b,x,a) \mapsto b.x.F(a)~.$$
As usual, $F$ induces a functor
$$F~:~\{\text{$A$-modules}\} \to \{\text{$B$-modules}\}~;~M \mapsto B\otimes_AM$$
which in turn induces morphisms of the algebraic $K$-groups
$$F~:~K_i(A) \to K_i(B)~~(i=0,1)~.$$

\bepr\label{Proposition4}
{\rm (i)} A ring endomorphism $F:A \to A$ determines an isomorphism
torsion structure
$$F~:~{\rm Iso}({\mathbb B}(A)) \to {\rm Iso}({\mathbb B}(A))$$
and hence a relative $F$-torsion
$$\begin{array}{l}
\tau^F~:~K^{iso}_1({\mathbb B}(A)) \to K_1(A)~;\\[1ex]
\hphantom{\tau^F~:~}
\tau^{iso}(f:M \to N) \mapsto \tau^F(f)~=~\tau(F(f)^{-1}f:M \to M)~.
\end{array}$$
{\rm (ii)} If $F=F^2:A \to A$ the isomorphism torsion structure
is idempotent and
$$\aligned
{}&F^2~=~F~:~K_1(A) \to K_1(A)~,\\
{}&{\rm im}(\tau^F:K^{iso}_1({\mathbb B}(A)) \to K_1(A))~
=~{\rm im}(1-F:K_1(A) \to K_1(A))\\
{}&\hphantom{{\rm im}(\tau^F:K^{iso}_1({\mathbb B}(A)) \to K_1(A))~}
=~{\rm ker}(F:K_1(A) \to K_1(A))~,\\
{}&K_1(A) ~=~{\rm im}(F:K_1(A) \to K_1(A))\oplus {\rm im}(\tau^F)~.
\endaligned$$
\hfill$\qs$
\enpr

\begin{defi}\label{Definition2}(Ranicki \cite{ranitor}, p.211)\\
The {\it isomorphism torsion} of a contractible finite chain complex $C$
in an additive category ${\mathbb A}$ is
$$\tau^{iso}(C)~=~\tau^{iso}(d+\Gamma:C_{odd} \to C_{even}) \in K^{iso}_1({\mathbb A})$$
with $d+\Gamma$ the isomorphism in ${\mathbb A}$
defined for any chain contraction $\Gamma:0 \simeq 1:C \to C$ by
$$d+\Gamma~=~\begin{pmatrix}
d & 0 & 0 & \dots \\
\Gamma & d & 0 & \dots \\
0 & \Gamma & d & \dots \\
\vdots & \vdots & \vdots \end{pmatrix}~:
C_{odd}~=~C_1 \oplus C_3 \oplus C_5 \oplus \dots \to
C_{even}~=~C_0 \oplus C_2 \oplus C_4 \oplus \dots~.$$
\hfill$\qs$
\end{defi}

A chain map $f:D \to E$ of finite chain complexes in ${\mathbb A}$ is a
chain equivalence if and only if the algebraic mapping cone $C(f)$
is contractible.

\begin{rema}\label{Remark2}
We refer to \S4 of \cite{ranitor} for the precise definition of the
isomorphism torsion of a chain equivalence $f:D \to E$
of finite chain complexes in ${\mathbb A}$, which involves the
{\it sign} pairing
$$\begin{array}{l}
\epsilon~:~K_0({\mathbb A}) \times K_0({\mathbb A}) \to K^{iso}_1({\mathbb A})~;\\[1ex]
([M],[N]) \mapsto \epsilon(M,N)~=~
\tau^{iso}(\begin{pmatrix} 0 & 1 \\ 1 & 0 \end{pmatrix}:M \oplus N \to N \oplus M)~.
\end{array}$$
The isomorphism torsion of a chain equivalence $f$ is of the form
$$\tau^{iso}(f)~=~\tau^{iso}(C(f))+\epsilon(M,N) \in K^{iso}_1({\mathbb A})
\eqno{(*)}$$
with $M,N$ objects in ${\mathbb A}$. The isomorphism torsion depends only
on the chain homotopy class of $f$. The sign term in $(*)$ is
necessary in order to ensure that the isomorphism torsion of the
composite of chain equivalences $f:D \to E$, $g:E \to F$ have the
logarithmic property
$$\tau^{iso}(gf)~=~\tau^{iso}(f)+\tau^{iso}(g) \in K^{iso}_1({\mathbb A})~.$$
For any isomorphism torsion structure $F$ on ${\mathbb A}$ and any objects
$M,N$ in $\mathbb A$
$$ F\begin{pmatrix} 0 & 1 \\ 1 & 0 \end{pmatrix}~=~
 \begin{pmatrix} 0 & 1 \\ 1 & 0 \end{pmatrix}~:~M \oplus N \to N \oplus M~.$$
Thus the $F$-torsion of a sign is
$$\tau^F(\epsilon(M,N))~=~0 \in K_1({\mathbb A})$$
and the sign term in $(*)$ can be ignored for $F$-torsion purposes.
\hfill$\qs$
\end{rema}

\section{The Main Theorem}
\label{s:main}

The terminology used to deal
with polynomial extensions is developed in \ref{su:defin}.  In sections
\ref{su:whpoly}, \ref{su:whpower}, \ref{su:whlaurent} we recall the
splitting theorems for the Whitehead group of a twisted polynomial
extension $A_{\rho}[z]$, the power series ring $A_{\rho}[[z]]$ and the
Laurent ring $A_{\rho}[z,z^{-1}]$, paving the way to the proof in
section \ref{su:whnovikov} of the Main Theorem on the Whitehead group
of the Novikov ring $A_{\rho}((z))$.

\subsection{Modules, matrices, polynomial extensions etc.}
\label{su:defin}

We shall work with left modules over rings.
Let $A$ be a ring with a unit (non-commutative in general).
An $A$-module morphism of finite direct sums of $A$-modules
$$f~:~M_1 \oplus M_2 \oplus \dots \oplus M_q \to
N_1 \oplus N_2 \oplus \dots \oplus N_p$$
can be identified with the $p \times q$ matrix
$f=(f_{ij})_{1 \leq i \leq p,1\leq j \leq q}$
with entries $A$-module morphisms $f_{ij}:M_j \to N_i$, such that
$$f(x_1,x_2,\dots,x_p)~=~
(\sum\limits^q_{j=1}f_{1j}(x_j),\sum\limits^q_{j=1}f_{2j}(x_j),\dots,
\sum\limits^q_{j=1}f_{pj}(x_j))~.$$
An $A$-module morphism $A \to A$ can be identified with an element of $A$,
using the isomorphism of additive groups
$$A \to {\rm Hom}_A(A,A)~;~ a \mapsto (x \mapsto xa)~.$$
Thus for $M_1=M_2=\dots=M_q=N_1=N_2=\dots=N_p=A$ we can identify a
morphism of f.g. free $A$-modules $f:A^q \to A^p$
with a $p \times q$ matrix $(f_{ij})$ with entries
$$f_{ij} \in {\rm Hom}_A(A,A)~=~A~.$$

We shall be working with the twisted polynomial ring
$A_\rho[z]$ and its extensions $A_{\rho}[[z]], A_{\rho}((z))$, which are defined for an
automorphism $\rho:A \to A$ of the ring $A$ and an indeterminate $z$ over $A$
such that $az=z\rho(a)$ for $a \in A$, as in the Introduction.

A {\it $\rho$-morphism} of $A$-modules $M,N$ is a morphism $f:M\to N$
of the underlying additive groups such that
$$f(ax)~=~\rho(a)f(x)~~(a\in A)~.$$
The composite of $\rho$-morphisms $f:M \to N$, $g:N \to P$ is a
$\rho^2$-morphism $gf:M \to P$.

For any $A$-module $M$ and $k \in \zz$ let $z^kM$ be the $A$-module
with elements $z^kx$ ($x\in M$) and
$$z^kx+z^ky~=~z^k(x+y)~,~a(z^kx)~=~z^k(\rho^k(a)x)~~(x,y \in M,a \in A)~.$$
The function
$$\zeta^k~:~M \to z^kM~;~x \mapsto z^kx$$
is a $\rho^{-k}$-isomorphism.

A $\rho^k$-morphism $f:M\to N$ of $A$-modules determines an $A$-module morphism
$$z^kf~:~M\to z^kN~;~x \mapsto z^kf(x)$$
and every $A$-module morphism $M \to z^kN$ is of this form.
Thus there is no essential difference between $\rho^k$-morphisms
$M\to N$ and $A$-module morphisms $M \to z^kN$.

For any f.g. free $A$-module $A^n$ there is defined a $\rho$-isomorphism
$$\begin{array}{l}
\theta_n~=~\rho\oplus \rho \dots \oplus \rho~:~A^n \to A^n ~;\\[1ex]
\hspace*{15mm} (a_1,a_2,\dots,a_n) \mapsto (\rho(a_1),\rho(a_2),\dots,\rho(a_n))
\end{array}$$
or equivalently an $A$-module isomorphism
$$z\theta_n ~:~A^n \to zA^n~;~(a_1,a_2,\dots,a_n) \mapsto z(\rho(a_1),\rho(a_2),\dots,\rho(a_n))~.$$

For any automorphism $\rho:A \to A$ and $k \in \zz$ the ring morphism
$\rho^k:A \to A$ induces a functor
$$(\rho^k)_!~:~\{\hbox{$A$-modules}\} \to \{\hbox{$A$-modules}\}~;~
M \mapsto (\rho^k)_!M$$
with a natural $A$-module isomorphism
$$(\rho^k)_!M \to z^kM~;~(a,x) \mapsto z^kax~.$$
For the inclusion $i:A \to A_{\rho}[z]$ we write the induced
$A_{\rho}[z]$-module as
$$A_{\rho}[z]\otimes_AM~=~M_{\rho}[z]~=~\sum\limits^{\infty}_{j=0}z^jM~.$$
For any $A$-module $M$
$$A_{\rho}[z,z^{-1}] \otimes_A M~=~M_{\rho}[z,z^{-1}]~=~
\sum\limits^{\infty}_{j=-\infty}z^jM$$
and for a f.g. $A$-module $M$
$$\begin{array}{l}
A_{\rho}[[z]]\otimes_A M~=~M_{\rho}[[z]]~=~\prod\limits^{\infty}_{j=0}z^jM~,\\
A_{\rho}((z)) \otimes_A M~=~M_{\rho}((z))~=~
\sum\limits^{-1}_{j=-\infty}z^jM \oplus \prod\limits^{\infty}_{k=0}z^kM~.
\end{array}$$
\begin{defi}\label{Definition5} Let $R$ be one of the rings
$$A_{\rho}[z]~~,~~A_{\rho}[[z]]~~,~~A_{\rho}[z,z^{-1}]~~,~~A_{\rho}((z))~.$$
(i) A f.g. projective $R$-module is {\it $A$-induced}
if it is of the form $R\otimes_AP$ for a f.g. projective $A$-module $P$.\\
(ii) Let ${\mathbb P}_A(R)$ be the additive category of $A$-induced f.g. projective
$R$-modules.$\qs$
\end{defi}

In particular, every f.g. free $R$-module is $A$-induced from a f.g. free
$A$-module
$$R^n~=~R\otimes_AA^n~.$$
If $P$ is a f.g. projective $A$-module and $Q$ is any $A$-module it is possible
to write every $A_{\rho}[z]$-module morphism $f:P_{\rho}[z] \to Q_{\rho}[z]$
as a polynomial
$$\begin{array}{ll}
f~=~\sum\limits^{\infty}_{j=0}z^jf_j~:&
P_{\rho}[z]~=~\sum\limits^{\infty}_{k=0}z^kP \to
Q_{\rho}[z]~=~\sum\limits^{\infty}_{\ell=0}z^{\ell}Q~;\\
&\sum\limits^{\infty}_{k=0}z^kx_k \mapsto
\sum\limits^{\infty}_{j=0}\sum\limits^{\infty}_{k=0}z^{j+k}f_j(x_k)
\end{array}$$
with coefficients $\rho^j$-morphisms $f_j :P \to Q$ such that
$\{j \geq 0\,\vert\,f_j \neq 0\}$ is finite. Similarly for
$A_{\rho}[[z]]$, $A_{\rho}[z,z^{-1}]$, $A_{\rho}((z))$.
\medskip

\bepr\label{Proposition5}
Let $R,{\mathbb P}_A(R)$ be as in Definition \ref{Definition5}.
The forgetful functor
$${\mathbb P}_A(R) \to {\mathbb P}(R)$$
from the additive category of $A$-induced f.g.  projective $R$-modules to
the additive category of f.g.  projective $R$-modules induces an
isomorphism of Whitehead groups
$$K_1({\mathbb P}_A(R))~\cong~K_1({\mathbb P}(R))~=~K_1(R)~.$$
\enpr
\noindent{\it Proof}. For any f.g. projective $A$-modules $P,Q$
$${\rm Hom}_{{\mathbb P}_A(R)}(R\otimes_AP,R\otimes_AQ)~=~
{\rm Hom}_R(R\otimes_AP,R\otimes_AQ)~.$$
\hfill$\qs$

The forgetful functor ${\mathbb P}_A(R) \to {\mathbb P}(R)$ is an equivalence
of additive categories, but the induced morphism
$K^{iso}_1({\mathbb P}_A(R)) \to K^{iso}_1({\mathbb P}(R))$ is not an
isomorphism (cf. Remark \ref{Remark1}).

\bepr\label{Proposition6} Let $R$ be one of the rings
$$A_{\rho}[z]~~,~~A_{\rho}[[z]]~.$$
{\rm (i)} The inclusion $A \to R$ is split by the augmentation map
$$R \to A~;~\sum\limits^{\infty}_{j=0}a_jz^j \mapsto a_0~.$$
{\rm (ii)} The composite of the augmentation
and the inclusion is an idempotent endomorphism
$$F~:~R \to R~;~\sum\limits^{\infty}_{j=0}a_jz^j \mapsto a_0~.$$
{\rm (iii)} The Whitehead group of $R$ splits as
$$K_1(R)~=~K_1(A) \oplus NK_1(R)$$
with
$$\begin{array}{ll}
NK_1(R)&=~{\rm ker}(K_1(R) \to K_1(A))\\[1ex]
&=~{\rm im}(1-F:K_1(R)\to K_1(R))\\[1ex]
&=~{\rm ker}(F:K_1(R)\to K_1(R))~.
\end{array}$$
{\rm (iv)} $F$ determines an idempotent isomorphism structure
$$F~:~{\rm Iso}({\mathbb P}_A(R)) \to {\rm Iso}({\mathbb P}_A(R))$$
with a relative $F$-torsion
$$\begin{array}{l}
\tau^F~:~K^{iso}_1({\mathbb P}_A(R)) \to NK_1(R)~;\\[1ex]
\tau^{iso}(f:M \to N) \mapsto \tau^F(f)=\tau(F(f)^{-1}f:M \to M)~.
\end{array}$$
\enpr
\noindent{\it Proof.} (i), (ii), (iii) Clear.\\
(iv) Every morphism in ${\mathbb P}_A(R)$ is of the form
$$f~=~\sum\limits^{\infty}_{j=0}f_jz^j~:~M~=~R\otimes_A P \to
N~=~R\otimes_A Q$$
for some f.g. projective $A$-modules $P,Q$ and $\rho^j$-morphisms
$f_j:P \to Q$. If $f$ is an isomorphism then so is
$$F(f)~=~f_0~:~M~=~R\otimes_A P \to N~=~R\otimes_A Q~.$$
The idempotent isomorphism structure is defined by
$$F~:~{\rm Iso}({\mathbb P}_A(R)) \to {\rm Iso}({\mathbb P}_A(R))~;~
\tau(f:M \to N) \mapsto \tau(f_0:M \to N)~.\eqno{\qs}$$

In the situation of Proposition \ref{Proposition6}
we shall write the natural isomorphism as
$$B~=~B_1 \oplus B_2~:~K_1(R) \to
K_1(A) \oplus NK_1(R)$$
and the inverse isomorphism as
$$B^{-1}~=~C~=~C_1 \oplus C_2~:~K_1(A) \oplus NK_1(R) \to
K_1(R)~.$$
The components
$$B_1~:~K_1(R)\to K_1(A) ~,~C_1~:~K_1(A) \to K_1(R)$$
are induced by  the augmentation $R \to A;z \mapsto 0$
and the inclusion $A \to R$, with
$$\begin{array}{l}
B_1C_1~=~1~:~K_1(A) \to K_1(A)~,\\[1ex]
C_1B_1~=~F~:~K_1(R) \to K_1(R)~.
\end{array}$$
The components
$$B_2~:~K_1(R)\to NK_1(R) ~,~C_2~:~NK_1(R) \to K_1(R)$$
are the natural surjection and injection, with
$$\begin{array}{l}
B_2C_2~=~1~:~NK_1(R) \to NK_1(R)~,\\[1ex]
C_2B_2~=~1-F~:~K_1(R) \to K_1(R)~,\\[1ex]
B_2C_1~=~0~:~K_1(A) \to NK_1(R)~,\\[1ex]
B_1C_2~=~0~:~NK_1(R) \to K_1(A)~.
\end{array}$$

\subsection{The Whitehead group of $A_{\rho}[z]$}
\label{su:whpoly}
By Proposition \ref{Proposition6} we have
$$K_1(A_{\rho}[z])~=~K_1(A) \oplus NK_1(A_{\rho}[z])$$
with
$$\begin{array}{ll}
NK_1(A_{\rho}[z])&=~{\rm ker}(K_1(A_{\rho}[z]) \to K_1(A))\\[1ex]
&=~{\rm im}(1-F:K_1(A_{\rho}[z])\to K_1(A_{\rho}[z]))\\[1ex]
&=~{\rm ker}(F:K_1(A_{\rho}[z])\to K_1(A_{\rho}[z]))
\end{array}$$
where
$$F~:~A_{\rho}[z] \to A_{\rho}[z]~;~
\sum\limits^{\infty}_{j=0}a_jz^j \mapsto a_0~.$$
We shall now recall the identification of $NK_1(A_{\rho}[z])$ with the
reduced $\rho$-nilpotent class group $\widetilde{\rm Nil}_0(A,\rho)$.

A $\rho$-endomorphism $\nu:P \to P$ is {\it nilpotent} if for some $k \geq 1$
$$\nu^k~=~0~:~P \to P~.$$

\bepr\label{Proposition7}
{\rm (i)} A linear morphism of $A$-induced f.g. projective $A_{\rho}[z]$-modules
$$\alpha~=~\alpha_0+z\alpha_1~:~Q_{\rho}[z] \to R_{\rho}[z]$$
is an isomorphism if and only if $\alpha_0:Q \to R$ is an $A$-module
isomorphism and $(\alpha_0)^{-1}\alpha_1:Q \to Q$ is a nilpotent
$\rho$-endomorphism.\\
{\rm (ii)} {\rm (}Higman linearization trick.{\rm)}
Every element of $K_1(A_{\rho}[z])$ is the torsion $\tau(\alpha)$
of a linear automorphism of a f.g. free $A_{\rho}[z]$-module.\hfill$\qs$
\enpr

The {\it nilpotent class group} $\nil$ is
the class group of the exact category ${\rm Nil}(A,\rho)$ with objects
$$(P=\hbox{\rm f.g. projective $A$-module},\nu:P \to P~\hbox{\rm
nilpotent $\rho$-endomorphism})~,$$
that is
$${\rm Nil}_0(A,\rho)~=~K_0({\rm Nil}(A,\rho))~.$$
The {\it reduced nilpotent class group} $\widetilde{\mbox{\rm Nil}}_0(A,\rho)$ is
defined to be
$$\widetilde{\mbox{\rm Nil}}_0(A,\rho)~=~{\rm coker}(K_0(A) \to
{\rm Nil}_0(A,\rho))$$
with
$$K_0(A) \to {\rm Nil}_0(A,\rho) ~;~ [P]\mapsto [P,0] ~,$$
and is such that
$${\rm Nil}_0(A,\rho)~=~\widetilde{\mbox{\rm Nil}}_0(A,\rho) \oplus K_0(A)~.$$
It is also possible to view $\widetilde{\mbox{\rm Nil}}_0(A,\rho)$ as
$$\widetilde{\mbox{\rm Nil}}_0(A,\rho)~=~{\rm coker}(K_0({\mathbb Z}) \to
K_0(\widetilde{\mbox{\rm Nil}}(A,\rho)))$$
with $\widetilde{\mbox{\rm Nil}}(A,\rho)$ the full subcategory of
${\rm Nil}(A,\rho)$ with objects $(F, \nu)$ such that $F$ is
a f.g.  free $A$-module, and
$$K_0({\mathbb Z})~=~{\mathbb Z} \to K_0(\widetilde{\mbox{\rm Nil}}(A,\rho))~;~
[{\mathbb Z}^n] \mapsto [A^n,0]~.$$

\begin{theo}\label{Theorem1}
{\rm (Bass, Heller and Swan \cite{bhs},
Bass \cite{bass}, Farrell and Hsiang \cite{farhsi})}\newline
The function
$$\widetilde{\rm Nil}_0(A,\rho) \to NK_1(A_{\rho}[z])~;~
[P,\nu] \mapsto \tau(1-z\nu:P_{\rho}[z] \to P_{\rho}[z])$$
is an isomorphism.\hfill$\qs$
\end{theo}

The following construction of nilpotent classes from isomorphisms
of $A$-induced f.g. projective $A_{\rho}[z]$-modules will be generalized
in section \ref{su:whpower} to a construction of Witt vector classes from isomorphisms
of $A$-induced f.g. projective $A_{\rho}[[z]]$-modules.

\bepr\label{Proposition8}
{\rm (i)} An isomorphism of $A$-induced f.g. projective $A_{\rho}[z]$-modules
$$\alpha~=~\sum\limits^k_{j=0}z^j\alpha_j~:~Q_{\rho}[z] \to R_{\rho}[z]$$
determines an object $(P,\nu)$ in ${\rm Nil}(A,\rho)$, with
$$\begin{array}{l}
P~=~{\rm coker}(\alpha\vert:z^{-k-1}Q_{\rho}[z^{-1}] \to z^{-1}R_{\rho}[z^{-1}])~,\\[1ex]
\nu~:~P \to P~;~[x] \mapsto [zx]~.
\end{array}$$
{\rm (ii)} For a linear isomorphism
$$\alpha~=~\alpha_0+z\alpha_1~:~Q_{\rho}[z] \to R_{\rho}[z]$$
the object in {\rm (i)} is given up to isomorphism by
$$(P,\nu)~=~(Q,-(\alpha_0)^{-1}\alpha_1)~.$$
{\rm (iii)} The torsion of an automorphism of an $A$-induced f.g. projective
$A_{\rho}[z]$-module
$$\alpha~=~\sum^k_{j=0}z^j\alpha_j~:~Q_{\rho}[z] \to Q_{\rho}[z]$$
splits as
$$\tau(\alpha)~=~\tau(\alpha_0)+\tau(1-z\nu:P_{\rho}[z] \to P_{\rho}[z])
\in K_1(A_{\rho}[z])$$
with $(P,\nu)$ as in {\rm (i)}.\\
{\rm (iv)} The isomorphism
$$B~=~B_1\oplus B_2~:~
K_1(A_{\rho}[z])\to K_1(A) \oplus \widetilde{\mbox{\rm Nil}}_0(A,\rho)$$
has components
$$\aligned
{}&B_1~:~K_1(A_{\rho}[z])\to K_1(A) ~;~
\tau(\sum\limits^{\infty}_{j=0}\alpha_jz^j:A_{\rho}[z]^n \to A_{\rho}[z]^n)
\mapsto \tau(\alpha_0:A^n \to A^n)~,\\[1ex]
{}&B_2~:~K_1(A_{\rho}[z])\to \widetilde{\mbox{\rm Nil}}_0(A,\rho) ~;~
\tau(\sum\limits^k_{j=0}z^j\alpha_j:A_{\rho}[z]^n \to A_{\rho}[z]^n)
\mapsto [P,\nu]
\endaligned$$
and the inverse isomorphism
$$B^{-1}~=~C~=~C_1\oplus C_2~:~K_1(A) \oplus
\widetilde{\mbox{\rm Nil}}_0(A,\rho) \to K_1(A_{\rho}[z])$$
has components
$$\aligned
{}&C_1~:~K_1(A) \to K_1(A_{\rho}[z])~;~\tau(\alpha_0:A^n \to A^n) \mapsto
\tau(\alpha_0:A_{\rho}[z]^n \to A_{\rho}[z]^n)~,\\
{}&C_2~:~\widetilde{\mbox{\rm Nil}}_0(A,\rho) \to K_1(A_{\rho}[z])~;~
[P,\nu] \mapsto \tau(1-z\nu:P_{\rho}[z] \to P_{\rho}[z])~.
\endaligned$$
{\rm (v)} For $F:A_{\rho}[z] \to A_{\rho}[z];z \mapsto z$ the relative
$F$-torsion function of \ref{Proposition6}
$$\tau^F~:~K^{iso}_1({\mathbb P}_A(A_{\rho}[z])) \to NK_1(A_{\rho}[z])~=~
\widetilde{\mbox{\rm Nil}}_0(A,\rho)$$
sends the isomorphism torsion of an isomorphism of $A$-induced f.g. projective
$A_{\rho}[z]$-modules
$$\alpha~=~\sum^k_{j=0}z^j\alpha_j~:~Q_{\rho}[z] \to R_{\rho}[z]$$
to the class of the nilpotent object $(P,\nu)$ in {\rm (i)}
$$\tau^F(\alpha)~=~[P,\nu] \in \widetilde{\mbox{\rm Nil}}_0(A,\rho)~.$$
\enpr
\noindent{\it Proof.} This is standard except for the isomorphism torsion
interpretation in (v), which is just an interpretation of
(i)--(iv) in terms of the material in section
\ref{s:class}. See Proposition \ref{Proposition13} for a detailed proof
that the $A$-module $P$ in (i) is f.g. projective, and that $\nu:P \to P$
is a nilpotent $\rho$-morphism.
\hfill$\qs$
\medskip

\begin{defi}\label{Definition6}
The {\it nilpotent torsion} is the group morphism given by the
relative $F$-torsion construction of Proposition \ref{Proposition8} (v)
$$\nu~=~\tau^F~:~K^{iso}_1({\mathbb P}_A(A_{\rho}[z]))
\to \widetilde{\mbox{\rm Nil}}_0(A,\rho)~;~
\tau^{iso}(\alpha) \mapsto [P,\nu]~.\eqno{\qs}$$
\end{defi}

The split surjection $B_2$ in Proposition \ref{Proposition8} (iv) is the
composite
$$B_2~:~K_1(A_{\rho}[z]) \to K^{iso}_1({\mathbb P}_A(A_{\rho}[z]))
\rTo^{\nu} \widetilde{\mbox{\rm Nil}}_0(A,\rho)~.$$

\begin{rema}\label{Remark3}
In view of \ref{Remark2}, it is also possible to define the nilpotent
torsion $\nu(f) \in \widetilde{\mbox{\rm Nil}}_0(A,\rho)$ of a chain
equivalence $f:D \to E$ of $A$-induced f.g. projective $A_{\rho}[z]$-module
chain complexes, which depends only on the chain homotopy class of $f$
and has the logarithmic property $\nu(fg)=\nu(f)+\nu(g)$.\hfill$\qs$
\end{rema}

\subsection{The Whitehead group of $A_{\rho}[[z]]$}
\label{su:whpower}
The splitting of $K_1(A_{\rho}[z])$ in section
\ref{su:whpoly} is now extended to $K_1(A_{\rho}[[z]])$. However,
there is an essential difference between the
decompositions
$$\begin{array}{l}
K_1(A_{\rho}[z])~=~K_1(A) \oplus NK_1(A_{\rho}[z])~,\\[1ex]
K_1(A_{\rho}[[z]])~=~K_1(A) \oplus NK_1(A_{\rho}[[z]])
\end{array}$$
in that there is no analogue for $K_1(A_{\rho}[[z]])$ of Higman linearization
for $K_1(A_{\rho}[z])$.

 As in Proposition \ref{Proposition6} we have
$$\begin{array}{ll}
NK_1(A_{\rho}[[z]])&=~{\rm ker}(K_1(A_{\rho}[[z]]) \to K_1(A))\\[1ex]
&=~{\rm im}(1-F:K_1(A_{\rho}[[z]])\to K_1(A_{\rho}[[z]]))\\[1ex]
&=~{\rm ker}(F:K_1(A_{\rho}[[z]])\to K_1(A_{\rho}[[z]]))
\end{array}$$
where
$$F~:~A_{\rho}[[z]] \to A_{\rho}[[z]]~;~
\sum\limits^{\infty}_{j=0}a_jz^j \mapsto a_0~.$$

\bepr\label{Proposition9}
An $A_{\rho}[[z]]$-module morphism of $A$-induced f.g.  projective
$A_{\rho}[[z]]$-modules
$$\alpha~=~\sum^{\infty}_{j=0}z^j\alpha_j~:~P_{\rho}[[z]] \to Q_{\rho}[[z]]$$
is an isomorphism if and only if $\alpha_0:P \to Q$ is an $A$-module
isomorphism.$\qs$
\enpr

\begin{defi} \label{Definition7} (i)
A {\it Witt vector} in $A_{\rho}[[z]]$ is a unit of the type
$$w~=~1 +\sum\limits^{\infty}_{j=1}a_jz^j \in A_{\rho}[[z]]^{\bullet}~.$$
(ii) Let $W_1(A,\rho)\subseteq K_1(A_{\rho}[[z]])$ be the subgroup
of the torsions $\tau(w)$ of Witt vectors $w$.
\hfill$\qs$
\end{defi}

\bepr\label{Proposition10} {\rm (Pajitnov \cite{patou})}\newline
$$NK_1(A_{\rho}[[z]])~=~ W_1(A,\rho) \subseteq
K_1(A_{\rho}[[z]]))$$
\enpr
\noindent{\it Proof.} By Proposition \ref{Proposition9} an endomorphism of
a f.g. free $A_{\rho}[[z]]$-module of rank $n$
$$\alpha~=~\sum\limits^{\infty}_{i=0}\alpha_iz^i~:~
A_{\rho}[[z]]^n \to A_{\rho}[[z]]^n$$
is an automorphism if and only if $\alpha_0:A^n \to A^n$ is an
$A$-module automorphism, in which case the torsion of $\alpha$ is the sum
$$\tau(\alpha)~=~\tau(\alpha_0)+\tau(\beta) \in K_1(A_{\rho}[[z]])$$
with
$$\beta~=~(\alpha_0)^{-1}\alpha~=~1+\sum\limits^{\infty}_{i=1}\beta_iz^i~:~
A_{\rho}[[z]]^n \to A_{\rho}[[z]]^n~.$$
The diagonal entries in the matrix $\beta=(\beta_{jk})$ are Witt vectors
$$\beta_{jj} \in 1+zA_{\rho}[[z]] \subset A_{\rho}[[z]]^{\bullet}~~(1 \leq j \leq n)$$
so that $\beta$ can be reduced by elementary row operations to an
upper triangular matrix with diagonal entries Witt vectors
$$\gamma~=~\begin{pmatrix}
w_1 & * & \dots & * \\
0 & w_2 & \dots & * \\
\vdots & \vdots & \ddots & \vdots \\
0 & 0 & \dots & w_n \end{pmatrix}~:~A_{\rho}[[z]]^n \to A_{\rho}[[z]]^n~.$$
Thus
$$\aligned
\tau(\alpha)~&=~\tau(\alpha_0) + \tau(\beta) \\
&=~\tau(\alpha_0) + \tau(\gamma) ~=~\tau(\alpha_0) + \sum\limits^n_{j=1}\tau(w_j)\\
&\in  K_1(A_{\rho}[[z]])~.
\endaligned$$
In this terminology
$$\begin{array}{l}
B_1(\tau(\alpha))~=~\tau(\alpha_0) \in K_1(A)~,\\[1ex]
B_2(\tau(\alpha))~=~\tau((\alpha_0)^{-1}\alpha)~=~
\sum\limits^n_{j=1}\tau(w_j)\in W_1(A,\rho)~.
\end{array}
$$
$\qs$

The proof of the Main Theorem will make use of the following
construction of elements in $W_1(A,\rho)$, analogous to the
construction \ref{Proposition8} of elements in
$\widetilde{\mbox{\rm Nil}}_0(A,\rho)$.

\begin{defi} \label{Definition8}
The {\it Witt vector torsion} is the group morphism given by
$$\begin{array}{l}
w~:~K^{iso}_1({\mathbb P}_A(A_{\rho}[[z]]))
\to W_1(A,\rho)~;~\tau^{iso}(\alpha) \mapsto w(\alpha)~=~\tau((\alpha_0)^{-1}\alpha)~.
\end{array}$$
Equivalently, this is the relative $F$-torsion of Proposition \ref{Proposition6}
$$\tau^F~:~K^{iso}_1({\mathbb P}_A(A_{\rho}[[z]])) \to
NK_1(A_{\rho}[[z]]))~=~W_1(A,\rho)$$
with $F:A_{\rho}[[z]] \to A_{\rho}[[z]];z \mapsto 0$.
\hfill$\qs$
\end{defi}
\medskip

\begin{rema}\label{Remark4}
In view of \ref{Remark2}, it is also possible to define the Witt vector
torsion $w(f) \in W_1(A,\rho)$ of a chain
equivalence $f:D \to E$ of $A$-induced f.g. projective $A_{\rho}[[z]]$-module
chain complexes, which depends only on the chain homotopy class of $f$
and has the logarithmic property $w(fg)=w(f)+w(g)$.\hfill$\qs$
\end{rema}

\subsection{The Whitehead group of $A_{\rho}[z,z^{-1}]$}
\label{su:whlaurent}

The aim of this section is to describe the decomposition of
$K_1(A_{\rho}[z,z^{-1}])$ in a way which will serve as a model for our
main result on the decomposition of $K_1(A_{\rho}((z)))$.
\medskip

\bepr\label{Proposition13}
{\rm (i)} An isomorphism of $A$-induced f.g. projective $A_{\rho}[z,z^{-1}]$-modules
$$\alpha~=~\sum\limits^{k_-}_{j=-k_+}z^j\alpha_j~:~
Q_{\rho}[z,z^{-1}] \to R_{\rho}[z,z^{-1}]$$
determines an object $(P_+,\nu_+)$ in ${\rm Nil}(A,\rho^{-1})$ and an
object $(P_-,\nu_-)$ in ${\rm Nil}(A,\rho)$, with
$$\begin{array}{l}
P_+~=~{\rm coker}(\alpha_+:z^{k_+}Q_{\rho}[z] \to R_{\rho}[z])~,\\[1ex]
\nu_+~:~P_+ \to P_+~;~[x] \mapsto [z^{-1}x]~,\\[1ex]
P_-~=~{\rm coker}(\alpha_-:z^{-k_--1}Q_{\rho}[z^{-1}] \to z^{-1}R_{\rho}[z^{-1}])~,\\[1ex]
\nu_-~:~P_- \to P_-~;~[x] \mapsto [zx]~,
\end{array}$$
with $\alpha_+,\alpha_-$ restrictions of $\alpha$.\\
{\rm (ii)} The constructions of {\rm (i)} define group morphisms
$$\begin{array}{l}
\nu_+~:~K^{iso}_1({\mathbb P}_A(A_{\rho}[z,z^{-1}]))
\to \widetilde{\mbox{\rm Nil}}_0(A,\rho^{-1})~;~\tau^{iso}(\alpha) \mapsto [P_+,\nu_+]~,\\[1ex]
\nu_-~:~K^{iso}_1({\mathbb P}_A(A_{\rho}[z,z^{-1}]))
\to \widetilde{\mbox{\rm Nil}}_0(A,\rho)~;~\tau^{iso}(\alpha) \mapsto
[P_-,\nu_-]~.
\end{array}$$
\enpr
\noindent{\it Proof.}
(i) It is clear that $P_+$ is a f.g. $A_{\rho}[z]$-module, and that there is
defined an $A$-induced f.g. projective $A_{\rho}[z]$-module resolution
$$\begin{diagram}
0& \rTo& z^{k_+}Q_{\rho}[z] & \rTo^{\alpha_+} &R_{\rho}[z] &\rTo^{\pi_+}& P_+&
\rTo& 0
\end{diagram}
$$
with
$$\pi_+~=~{\rm projection}~:~R_{\rho}[z] \to P_+~,$$
and similarly for the $A_{\rho}[z^{-1}]$-module $P_-$ with a resolution
$$\begin{diagram}
0& \rTo& z^{-k_--1}Q_{\rho}[z^{-1}] & \rTo^{\alpha_-} &
z^{-1}R_{\rho}[z^{-1}] & \rTo^{\pi_-}& P_-& \rTo& 0~.
\end{diagram}$$
It follows from the commutative diagram
\begin{diagram}
0 & \rTo & z^{k_+}Q_{\rho}[z]\oplus z^{-k_--1}Q_{\rho}[z^{-1}]
& \rTo    & Q_{\rho}[z,z^{-1}]  & \rTo  & \sum\limits^{k_+-1}_{j=-k_-}z^jQ &\rTo & 0 \\
  &      &  \dTo^{\rm id}   &              & \dTo^{\alpha} &     & \dTo   &   \\
0 & \rTo & z^{k_+}Q_{\rho}[z]\oplus z^{-k_--1}Q_{\rho}[z^{-1}] &
\rTo^{\alpha_+\oplus \alpha_-}  &R_{\rho}[z,z^{-1}]&
\rTo& P_+ \oplus P_-  & \rTo  & 0
\end{diagram}
\noindent that there is defined an $A$-module isomorphism
$$\sum\limits^{k_+-1}_{j=-k_-}z^jQ~\cong~P_+ \oplus P_-~,$$
so that $P_+,P_-$ are f.g. projective $A$-modules.
The $\rho^{-1}$-endomorphism
$$z~:~R_{\rho}[z] \to R_{\rho}[z]~;~ x \mapsto zx$$
fits into the commutative diagram
\begin{diagram}
0 & \rTo & z^{k_+}Q_{\rho}[z] & \rTo^{\alpha_+}    & R_{\rho}[z]  & \rTo^{\pi_+}  & P_+ &\rTo & 0 \\
  &      &  \dTo^{z}   &              & \dTo^{z} &     & \dTo^{\nu_+}   &   \\
0 & \rTo & z^{k_+}Q_{\rho}[z] & \rTo^{\alpha_+}    & R_{\rho}[z]  & \rTo^{\pi_+}  & P_+ &\rTo & 0
\end{diagram}
Write the inverse of $\alpha$ as
$$\alpha^{-1}~=~\beta~=~\sum\limits^{\ell_-}_{j=-\ell_+}z^j\beta_j~:~
R_{\rho}[z,z^{-1}] \to Q_{\rho}[z,z^{-1}]$$
for some $\ell_+,\ell_- \geq 0$. For any $x \in R_{\rho}[z]$
$$z^{k_++\ell_+}x~=~\alpha\beta(z^{k_++\ell_+}x)~=~
\alpha_+(z^{k_+}\beta(z^{\ell_+}x)) \in
{\rm im}(\alpha_+:z^{k_+}Q_{\rho}[z] \to R_{\rho}[z])$$
so that
$$(\nu_+)^{k_++\ell_+}\pi_+(x)~=~0 \in P_+~.$$
Thus $\nu_+:P_+ \to P_+$ is nilpotent, with
$$(\nu_+)^{k_++\ell_+}~=~0~:~P_+ \to P_+~.$$
Similarly, the $\rho$-endomorphism
$$z^{-1}~:~z^{-1}R_{\rho}[z^{-1}] \to z^{-1}R_{\rho}[z^{-1}]~;~ x \mapsto z^{-1}x$$
fits into the commutative diagram
\begin{diagram}
0 & \rTo & z^{-k_--1}Q_{\rho}[z^{-1}] & \rTo^{\alpha_-}    &
z^{-1}R_{\rho}[z^{-1}]  & \rTo^{\pi_-}  & P_- &\rTo & 0 \\
  &      &  \dTo^{z^{-1}}   &              & \dTo^{z^{-1}} &     & \dTo^{\nu_-}   &   \\
0 & \rTo & z^{-k_--1}Q_{\rho}[z^{-1}] & \rTo^{\alpha_-}    &
z^{-1}R_{\rho}[z^{-1}]  & \rTo^{\pi_-}  & P_- &\rTo & 0
\end{diagram}
and
$$(\nu_-)^{k_-+\ell_-}~=~0~:~P_- \to P_-~.$$
(ii) The composite of isomorphisms of $A$-induced f.g. projective $A_{\rho}[z]$-modules
$$\begin{array}{l}
\alpha~=~\sum\limits^{k_+}_{j=-k_-}z^j\alpha_j~:~
Q_{\rho}[z,z^{-1}] \to R_{\rho}[z,z^{-1}]~,\\[1ex]
\alpha'~=~\sum\limits^{k'_+}_{j=-k'_-}z^j\alpha'_j~:~
R_{\rho}[z,z^{-1}] \to S_{\rho}[z,z^{-1}]
\end{array}$$
is an isomorphism of $A$-induced f.g. projective $A_{\rho}[z,z^{-1}]$-modules
$$\alpha''~=~\alpha'\alpha~=~\sum\limits^{k_++k_+'}_{j=-k_--k'_-}z^j\alpha''_j~:~
Q_{\rho}[z,z^{-1}] \to S_{\rho}[z,z^{-1}]$$
such that the corresponding nilpotent objects fit into exact sequences
$$\begin{array}{l}
0 \to (P_+,\nu_+) \to (P_+'',\nu_+'') \to (P_+',\nu_+') \to 0~,\\[1ex]
0 \to (P_-,\nu_-) \to (P_-'',\nu_-'') \to (P_-',\nu_-') \to 0~.
\end{array}$$
$\qs$

\begin{defi}\label{Definition10}
(i) The {\it automorphism class group} of $A,\rho$ is the class group
$${\rm Aut}_0(A,\rho)~=~K_0({\rm Aut}(A,\rho))$$
of the exact category ${\rm Aut}(A,\rho)$
of pairs $(P, \phi)$ with $P$ a f.g. projective $A$-module and $\phi: P\to P$
a $\rho$-isomorphism.\\
(ii) (Siebenmann \cite{sieben}) The {\it class-torsion group}
$K_1(A,\rho)$ is the quotient of ${\rm Aut}_0(A,\rho)$ by the subgroup
generated by the differences $(P,\phi)-(P',\phi')$ for which there
exists an isomorphism $h:P \to P'$ such that
$$\tau(h^{-1}{\phi'}^{-1}(zh)\phi:P \to zP \to zP' \to P' \to P)=0 \in
K_1(A)~.$$
\hfill$\qs$ \end{defi}

The group $K_1(A,\rho)$
fits into the long exact sequence
$$K_1(A)\rTo^{1-\rho} K_1(A)\rTo^i K_1(A,\rho) \rTo^j K_0(A)
\rTo^{1-\rho} K_0(A)$$
with
$$\begin{array}{l}
i~:~ K_1(A) \to K_1(A,\rho)~;~\tau(\alpha:A^n \to A^n) \mapsto
[A^n,\theta_n\alpha] - [A^n,\theta_n]~,\\[1ex]
j~:~K_1(A,\rho) \to K_0(A)~;~ [P,\phi] \mapsto [P]~.
\end{array}$$
\smallskip

\begin{theo}\label{Theorem3}
{\rm (Bass \cite{bass}, Farrell and Hsiang \cite{farhsi}, Siebenmann \cite{sieben})}\newline
There is a natural decomposition
$$K_1(A_{\rho}[z,z^{-1}])~=~K_1(A,\rho) \oplus \widetilde{\mbox{\rm Nil}}_0(A,\rho) \oplus
\widetilde{\mbox{\rm Nil}}_0(A,\rho^{-1})$$
with the map
$$C~=~C_1\oplus C_2\oplus C_3~:~
K_1(A,\rho) \oplus \widetilde{\mbox{\rm Nil}}_0(A,\rho) \oplus
\widetilde{\mbox{\rm Nil}}_0(A,\rho^{-1}) \to K_1(A_{\rho}[z,z^{-1}])$$
defined by
$$\begin{array}{l}
C_1~:~K_1(A,\rho)\to K_1(A_{\rho}[z,z^{-1}])~;\\[1ex]
\hspace*{15mm} (P,\phi)\mapsto
\tau(z\phi:P_\rho[z,z^{-1}]\to P_\rho[z,z^{-1}])~,\\[1ex]
C_2~:~\widetilde{\mbox{\rm Nil}}_0(A,\rho)\to K_1(A_{\rho}[z,z^{-1}])~;\\[1ex]
\hspace*{15mm} (P_+,\nu_+,)\mapsto
\tau(1-z\nu_+:(P_+)_\rho[z,z^{-1}]\to (P_+)_\rho[z,z^{-1}])~,\\[1ex]
C_3~:~\widetilde{\mbox{\rm Nil}}_0(A,\rho^{-1})\to K_1(A_{\rho}[z,z^{-1}])~;
\\[1ex]
\hspace*{15mm} (P_-,\nu_-)\mapsto
\tau(1-z^{-1}\nu_-:(P_-)_\rho[z,z^{-1}]\to (P_-)_\rho[z,z^{-1}])
\end{array}$$
an isomorphism.$\qs$
\end{theo}

The inverse isomorphism
$$C^{-1}~=~B~=~B_1\oplus B_2\oplus B_3~:~K_1(A_{\rho}[z,z^{-1}]) \to
K_1(A,\rho) \oplus \widetilde{\mbox{\rm Nil}}_0(A,\rho) \oplus
\widetilde{\mbox{\rm Nil}}_0(A,\rho^{-1})$$
is constructed as follows.

\bele\label{Lemma1}
{\rm (i)} The following sequence is an $A$-induced f.g.  projective
$A_{\rho}[z]$-module resolution of $P_+$\, :
\begin{diagram}
0&\rTo&(zP_+)_\rho [z]  &\rTo^{z\zeta^{-1}-\nu_+}& (P_+)_{\rho}[z] & \rTo^{\pi_+}& P_+ &\rTo& 0
\end{diagram}
with
$$\begin{array}{l}
\zeta^{-1}~:~zP_+ \to P_+~;~zx \mapsto x~~(\hbox{\it a $\rho$-isomorphism})~,\\
z\zeta^{-1}-\nu_+~:~(zP_+)_{\rho}[z] \to (P_+)_{\rho}[z] ~;~
\sum\limits^{\infty}_{j=1}z^jx_j \mapsto
\sum\limits^{\infty}_{j=1}z^jx_j-\sum\limits^{\infty}_{j=0}z^j\nu_+(x_{j+1})~,\\
\pi_+~:~(P_+)_{\rho}[z] \to P_+~;~
\sum\limits^{\infty}_{j=0}z^jx_j \mapsto \sum\limits^{\infty}_{j=0}(\nu_+)^j(x_j)~.
\end{array}
$$
{\rm (ii)} The two $A_{\rho}[z]$-module resolutions for $P_+$ are related by a
commutative diagram
\begin{diagram}
0&\rTo&z^{k_+}A_{\rho}[z]^n&\rTo^{\alpha_+}& A_{\rho}[z]^n &\rTo^{\pi_+} &P_+ &\rTo& 0 \\
 &    &\dTo^{g_+}&            & \dTo^{f_+} &     &\dTo^{\rm id}& & \\
0&\rTo&(zP_+)_\rho [z]  &\rTo^{z\zeta^{-1}-\nu_+}& (P_+)_{\rho}[z] & \rTo^{\widetilde \pi_+}& P_+ &\rTo& 0 \\
\end{diagram}
\noindent with
$$\begin{array}{l}
\widetilde \pi_+~:~(P_+)_{\rho}[z] \to P_+~;~\sum\limits^{\infty}_{j=0}z^jx_j \mapsto
\sum\limits^{\infty}_{j=0}(\nu_+)^j(x_j)~,\\[1ex]
f_+~:~A_{\rho}[z]^n \to (P_+)_{\rho}[z]~;~\sum\limits^{\infty}_{j=0}z^ja_j \mapsto
\sum\limits^{\infty}_{j=0}z^j\pi_+(a_j)~,\\[1ex]
g_+~:~z^{k_+}A_{\rho}[z]^n \to (zP_+)_{\rho}[z]~;~
\sum\limits^{\infty}_{j=k_+}z^jb_j \mapsto
\sum\limits^{\infty}_{j=1}\sum\limits^j_{i=1}z^i(\nu_+)^{j-i}(y_j)\\[1ex]
\hspace*{20mm}
(f_+\alpha_+(\sum\limits^{\infty}_{j=k_+}z^jb_j)=\sum\limits^{\infty}_{j=0}z^jy_j)~.
\end{array}$$
$\qs$
\enle

Regarding the commutative diagram in \ref{Lemma1} (ii) as a chain equivalence
of 1-dimensional $A$-induced f.g. projective $A_{\rho}[z]$-module chain complexes
$$(f_+,g_+)~:~C(\alpha_+:z^{k_+}A_{\rho}[z]^n \to A_{\rho}[z]^n)
\to C(z\zeta^{-1}-\nu_+:(zP_+)_\rho [z] \to (P_+)_{\rho}[z])$$
we have that the algebraic mapping cone is a short exact sequence of
$A$-induced f.g. projective $A_{\rho}[z]$-modules
$$C(f_+,g_+)~:~0 \to (z^{k_+}A^n)_{\rho}[z] \to (zP_+\oplus A^n)_{\rho}[z] \to (P_+)_{\rho}[z] \to 0~.$$
Choosing a splitting there is obtained an $A_{\rho}[z]$-module isomorphism
$$h~=~\sum\limits^{\infty}_{j=0}z^jh_j~:~
(P_+\oplus z^{k_+}A^n)_{\rho}[z]\rTo^{\cong}(zP_+ \oplus A^n)_{\rho}[z]~.$$
The components are $\rho^j$-morphisms
$$h_j~:~P_+ \oplus z^{k_+} A^n \to zP_+ \oplus A^n~~(j \geq 0)$$
with $h_0:P_+ \oplus z^{k_+} A^n \to zP_+ \oplus A^n$ an $A$-module isomorphism.
Use the $\rho$-isomorphism
$$z^{-1}h_0~:~P_+ \oplus z^{k_+}A^n \to P_+ \oplus z^{-1}A^n~;~x \mapsto z^{-1}(h_0(x))$$
and the $A$-module isomorphisms
$$\begin{array}{l}
z^{k_+}(\theta_n)^{k_+}~:~A^n \to z^{k_+} A^n~;~ (a_1,a_2,\dots,a_n)
\mapsto z^{k_+}(\rho^{k_+}(a_1),\rho^{k_+}(a_2),\dots,\rho^{k_+}(a_n))\\[1ex]
z\theta_n~:~z^{-1}A^n \to A^n~;~
z^{-1}(a_1,a_2,\dots,a_n)\mapsto (\rho(a_1),\rho(a_2),\dots,\rho(a_n))
\end{array}$$
to define a $\rho$-isomorphism
$$\phi~=~(1 \oplus z\theta_n)(z^{-1}h_0)
(1 \oplus z^{k_+}(\theta_n)^k)~:~P_+ \oplus A^n \to P_+ \oplus A^n~.$$

The components $B_i$ of $C=B^{-1}$ are given by
$$\begin{array}{l}
B_1~:~K_1(A_{\rho}[z,z^{-1}])\to K_1(A,\rho) ~;~ \tau(\alpha)\mapsto
[P_+\oplus A^n,\phi] - [A^n,\theta_n]~,\\[1ex]
B_2~:~K_1(A_{\rho}[z,z^{-1}])\to \widetilde{\mbox{\rm Nil}}_0(A,\rho) ~;~\tau(\alpha) \mapsto [P_-,\nu_-]~,\\[1ex]
B_3~:~K_1(A_{\rho}[z,z^{-1}])\to \widetilde{\mbox{\rm Nil}}_0(A,\rho^{-1})~;~\tau(\alpha) \mapsto [P_+,\nu_+]~.
\end{array}$$

\subsection{The Whitehead group of $A_{\rho}((z))$}
\label{su:whnovikov}

The splitting of $K_1(A_{\rho}[z,z^{-1}])$ in section
\ref{su:whlaurent} is now extended to $K_1(A_{\rho}((z)))$. However, as
for $K_1(A_{\rho}[z])$, $K_1(A_{\rho}[[z]])$ there is an essential
difference between the proofs of these results\,: there is no analogue
for $K_1(A_{\rho}((z)))$ of the Higman linearization trick by which
every element of $K_1(A_{\rho}[z,z^{-1}])$ can be represented by a
linear automorphism
$$\alpha~=~\alpha_0 + z\alpha_1~:~A_{\rho}[z,z^{-1}]^n \to A_{\rho}[z,z^{-1}]^n~.$$

We shall show that the map
$$\widehat C~=~\widehat C_1\oplus \widehat C_2\oplus \widehat C_3~:~
K_1(A,\rho) \oplus W_1(A,\rho) \oplus \widetilde{\mbox{\rm Nil}}_0(A,\rho^{-1})
\to K_1(A_{\rho}((z)))$$
defined by
$$\begin{array}{l}
\widehat C_1~:~K_1(A,\rho)\to K_1(A_{\rho}((z)))~;~ (P,\phi)\mapsto
\tau(z\phi:P_\rho((z))\to P_\rho((z)))~, \\[1ex]
\widehat C_2~:~W_1(A,\rho)\to K_1(A_{\rho}((z)))~;~ w\mapsto w~, \\[1ex]
\widehat C_3~:~\widetilde{\mbox{\rm Nil}}_0(A,\rho^{-1})\to K_1(A_{\rho}((z)))~;~ (P,\nu)\mapsto
\tau(1-z^{-1}\nu:P_\rho((z))\to P_\rho((z)))
\end{array}$$
is an isomorphism by constructing an explicit inverse
$$\widehat C^{-1}~=~\widehat B~=~\widehat B_1\oplus \widehat B_2\oplus \widehat B_3~:~
K_1(A_{\rho}((z))) \to
K_1(A,\rho) \oplus W_1(A,\rho) \oplus \widetilde{\mbox{\rm Nil}}_0(A,\rho^{-1})~.$$
The definition of $\widehat B$ is based on the constructions of
several auxiliary objects, to which we now proceed.

An element $\alpha\in K_1(A_{\rho}((z)))$ is represented by
an automorphism of a f.g. free $A_{\rho}((z))$-module
$$\alpha~=~\sum\limits^{\infty}_{j=-k}z^j\alpha_j~:~
A_{\rho}((z))^n \to A_{\rho}((z))^n$$
with coefficients $\rho^j$-morphisms $\alpha_j:A^n \to A^n$,
for some $k \geq 0$.
Write the inverse of $\alpha$ as
$$\alpha^{-1}~=~\beta~=~\sum\limits^{\infty}_{j=-\ell}z^j\beta_j~:~
A_{\rho}((z))^n \to A_{\rho}((z))^n$$
for some $\ell \geq 0$.

Now we can define the object which plays the most important role in
our construction\,: the $A$-module
$$P~=~{\rm coker}(\widetilde\alpha:z^kA_{\rho}[[z]]^n \to A_{\rho}[[z]]^n)$$
with $\widetilde\alpha$ the restriction of $\alpha$.

\bele\label{Lemma2}
{\rm (i)} $P$ is a f.g. projective $A$-module.\\
{\rm (ii)} The $\rho^{-1}$-endomorphism
$$z~:~P \to P~;~ x \mapsto zx$$
is nilpotent.
\enle
\noindent{\it Proof.}
(i) It is clear that $P$ is a f.g.  $A_{\rho}[[z]]$-module, and that there is
defined a f.g. free $A_{\rho}[[z]]$-module resolution
$$\begin{diagram}
0& \rTo& z^kA_{\rho}[[z]]^n & \rTo^{\widetilde\alpha} &A_{\rho}[[z]]^n &\rTo^{\pi}& P&
\rTo& 0
\end{diagram}
$$
with
$$\pi~=~{\rm projection}~:~A_{\rho}[[z]]^n \to P~.$$
To show that $P$ is a f.g.  projective $A$-module, consider the commutative diagram
\begin{diagram}
0 & \rTo & z^kA_{\rho}[[z]]^n & \rTo^{\widetilde\alpha}    & A_{\rho}[[z]]^n  & \rTo^{\pi}  & P &\rTo & 0 \\
  &      &  \dTo^{\rm id}   &              & \dTo^{\widetilde\beta} &     & \dTo   &   \\
0 & \rTo & z^kA_{\rho}[[z]]^n & \rTo       & z^{-\ell}A_{\rho}[[z]]^n  & \rTo  &
\sum_{j=-\ell}^{k-1}z^j A^n
  &\rTo & 0 \\
  &      &  \dTo^{\rm id}   &         & \dTo^{\gamma}  &     & \dTo   &   \\
0 & \rTo & z^kA_{\rho}[[z]]^n & \rTo^{\widetilde\alpha}    & A_{\rho}[[z]]^n  & \rTo^{\pi}  & P &\rTo & 0
\end{diagram}
\noindent with $\widetilde\beta$ the restriction of $\beta$, and $\gamma$ the
$A$-module morphism defined by
$$\begin{array}{c}
\gamma~:~z^{-\ell}A_{\rho}[[z]]^n \to A_{\rho}[[z]]^n~;~\sum\limits^{\infty}_{j=-\ell}z^jx_j
\mapsto \sum\limits^{\infty}_{j=0}z^jy_j\\[1ex]
(\alpha(\sum\limits^{\infty}_{j=-\ell}z^jx_j)~=~
\sum\limits^{\infty}_{j=-k-\ell}z^jy_j)~.
\end{array}$$
It follows from the identity
$$\gamma\widetilde \beta~=~(\alpha\beta)\vert_{A_{\rho}[[z]]^n}~=~1~:~
A_{\rho}[[z]]^n \to A_{\rho}[[z]]^n$$
that $P$ is a direct summand of the f.g. free $A$-module
$\sum\limits_{j=-\ell}^{k-1}z^j A^n$, and hence a f.g. projective $A$-module.
In fact, the $A$-module defined by
$$Q~=~{\rm coker}(\widetilde\beta:A_{\rho}[[z]]^n \to z^{-\ell}A_{\rho}[[z]]^n)$$
is such that there exists an $A$-module isomorphism
$$P \oplus Q~\cong~\sum_{j=-\ell}^{k-1}z^j A^n~.$$
(ii) Recall that $P$ is not only a left $A$-module, but also
a left $A_{\rho}[[z]]$-module. The left multiplication by the element
$z\in A_{\rho}[[z]]$ defines a $\rho^{-1}$-endomorphism of the $A$-module $A_{\rho}[[z]]^n$
$$z~:~A_{\rho}[[z]]^n \to A_{\rho}[[z]]^n~;~ x \mapsto zx$$
which induces a $\rho^{-1}$-endomorphism of the $A$-module $P$
$$\nu~:~P \to P~;~ x \mapsto zx$$
with a commutative diagram
\begin{diagram}
0 & \rTo & z^kA_{\rho}[[z]]^n & \rTo^{\widetilde\alpha}    & A_{\rho}[[z]]^n  & \rTo^{\pi}  & P &\rTo & 0 \\
  &      &  \dTo^{z}   &              & \dTo^{z} &     & \dTo^{\nu}   &   \\
0 & \rTo & z^kA_{\rho}[[z]]^n & \rTo^{\widetilde\alpha}    & A_{\rho}[[z]]^n  & \rTo^{\pi}  & P &\rTo & 0
\end{diagram}

\noindent
For any $x \in A_{\rho}[[z]]^n$
$$z^{k+\ell}x~=~\alpha\beta(z^{k+\ell}x)~=~\widetilde\alpha(z^k\beta(z^{\ell}x)) \in
{\rm im}(\widetilde\alpha:z^kA_{\rho}[[z]]^n \to A_{\rho}[[z]]^n)$$
so that
$$\nu^{k+\ell}\pi(x)~=~0 \in P$$
and
$$\nu^{k+\ell}~=~0~:~P \to P~.\eqno{\qs}$$

We thus obtain an element
$$(P,\nu)\in \widetilde{\mbox{\rm Nil}}_0(A,\rho^{-1})~.$$

Define an $A_{\rho}[[z]]$-module morphism
$$\pi~:~P_{\rho}[[z]]\to P~;~\sum\limits^{\infty}_{j=0}z^jx_j \mapsto
\sum\limits^{\infty}_{j=0}\nu^j(x_j)~.$$
Note that the right hand side of the formula is well defined
since $\nu$ is nilpotent, and therefore the sum contains
only a finite number of terms.

Identify
$$zP_\rho [[z]]~=~\prod\limits^{\infty}_{j=1}z^jP~~,~~
P_\rho [[z]]~=~\prod\limits^{\infty}_{j=0}z^jP~.$$

By analogy with \ref{Lemma1}\,:

\bele\label{Lemma3}
{\rm (i)} The following sequence is an $A$-induced
f.g. projective $A_{\rho}[[z]]$-module resolution of $P$\, :
\begin{diagram}
0&\rTo&zP_\rho [[z]]  &\rTo^{z\zeta^{-1}-\nu}& P_{\rho}[[z]] & \rTo^{\widetilde \pi}& P &\rTo& 0
\end{diagram}
with
$$\begin{array}{l}
\zeta^{-1}~:~zP \to P ~;~zx \mapsto x~,\\
z\zeta^{-1}-\nu~:~zP_{\rho}[[z]] \to P_{\rho}[[z]] ~;~
\sum\limits^{\infty}_{j=1}z^jx_j \mapsto
\sum\limits^{\infty}_{j=1}z^jx_j-\sum\limits^{\infty}_{j=0}z^j\nu(x_{j+1})~,\\
\widetilde \pi~:~P_{\rho}[[z]] \to P~;~
\sum\limits^{\infty}_{j=0}z^jx_j \mapsto \sum\limits^{\infty}_{j=0}\nu^j(x_j)~.
\end{array}
$$
{\rm (ii)} The two $A_{\rho}[[z]]$-module resolutions for $P$ are related by a
commutative diagram
\begin{diagram}
0&\rTo&z^kA_{\rho}[[z]]^n&\rTo^{\widetilde \alpha}& A_{\rho}[[z]]^n &\rTo^{\pi} &P &\rTo& 0 \\
 &    &\dTo^g&            & \dTo^f &     &\dTo^{\rm id}& & \\
0&\rTo&zP_\rho [[z]]  &\rTo^{z\zeta^{-1}-\nu}& P_{\rho}[[z]] & \rTo^{\widetilde \pi}& P &\rTo& 0 \\
\end{diagram}
\noindent with
$$\begin{array}{l}
f~:~A_{\rho}[[z]]^n \to P_{\rho}[[z]]~;~\sum\limits^{\infty}_{j=0}z^ja_j \mapsto
\sum\limits^{\infty}_{j=0}z^j\pi(a_j)~,\\[1ex]
g~:~z^kA_{\rho}[[z]]^n \to zP_{\rho}[[z]]~;~
\sum\limits^{\infty}_{j=k}z^jb_j \mapsto
\sum\limits^{\infty}_{j=1}\sum\limits^j_{i=1}z^i\nu^{j-i}(y_j)\\[1ex]
\hspace*{20mm}
(f\widetilde\alpha(\sum\limits^{\infty}_{j=k}z^jb_j)=\sum\limits^{\infty}_{j=0}z^jy_j)~.
\end{array}$$
\enle
\noindent{\it Proof.} (i) The sequence is part of a direct sum system
of $A$-modules
\begin{diagram}
zP_{\rho}[[z]]& \pile{\rTo^{z\zeta^{-1}-\nu} \\ \lTo_{\tau}} &
P_{\rho}[[z]]& \pile{\rTo^{\widetilde\pi} \\ \lTo_{\sigma}} & P~,
\end{diagram}
with the $A$-module morphisms defined by
$$\begin{array}{l}
\sigma~:~P \to P_{\rho}[[z]]~;~x \mapsto x~,\\
\tau~:~P_{\rho}[[z]] \to zP_\rho [[z]]~;~
\sum\limits^{\infty}_{j=0}z^jx_j \mapsto \sum\limits^{\infty}_{j=0}
z^{j+1}(\sum\limits^{\infty}_{k=0}\nu^k(x_{j+k+1}))
\end{array}$$
such that
$$\begin{array}{l}
\widetilde \pi \circ \sigma~=~{\rm id}~:~P \to P~,\\[1ex]
\tau \circ (z\zeta^{-1}-\nu)~=~{\rm id}~:~zP_{\rho}[[z]] \to zP_{\rho}[[z]]~,\\[1ex]
\sigma \circ \widetilde \pi + (z\zeta^{-1}-\nu) \circ \tau~=~{\rm id}~:~P_{\rho}[[z]] \to P_{\rho}[[z]]~.
\end{array}$$
(ii) By construction.$\qs$

Regard the commutative diagram of Lemma \ref{Lemma3} (ii)
as a chain equivalence of
1-dimensional $A$-induced f.g. projective $A_{\rho}[[z]]$-module chain complexes
$$(f,g,k)~:~C(\widetilde\alpha:z^kA_{\rho}[[z]]^n \to A_{\rho}[[z]]^n)
\to C(z\zeta^{-1}-\nu:zP_\rho[[z]] \to P_{\rho}[[z]])~.$$
By Definition \ref{Definition8} and Remark \ref{Remark4} we now have
a Witt vector torsion
$$w(f,g,k) \in W_1(A,\rho)~.$$
More precisely, the algebraic mapping cone is a short exact sequence of
$A$-induced f.g. projective $A_{\rho}[[z]]$-modules
$$C(f,g,k)~:~0 \to (z^kA^n)_{\rho}[[z]] \to (zP\oplus A^n)_{\rho}[[z]]
\to P_{\rho}[[z]] \to 0~.$$
Choosing a splitting there is obtained an $A_{\rho}[[z]]$-module isomorphism
$$h~=~\sum\limits^{\infty}_{j=0}z^jh_j~:~
(P\oplus z^kA^n)_{\rho}[[z]]\rTo^{\cong}(zP \oplus A^n)_{\rho}[[z]]$$
such that the components are $\rho^j$-morphisms
$$h_j~:~P \oplus z^k A^n \to zP \oplus A^n~~(j \geq 0)$$
with $h_0:P \oplus z^k A^n \to zP \oplus A^n$ an $A$-module isomorphism,
and
$$(h_0)^{-1}h~=~1+\sum\limits^{\infty}_{j=1}z^j(h_0)^{-1}h_j~
:~(P\oplus z^kA^n)_{\rho}[[z]]\to (P\oplus z^kA^n)_{\rho}[[z]]$$
an automorphism of an $A$-induced f.g. projective $A_{\rho}[[z]]$-module,
such that
$$\begin{array}{l}
\tau^{iso}(f,g,k)~=~\tau^{iso}(h)
\in K^{iso}_1({\mathbb P}(A_{\rho}[[z]]))\\[1ex]
w(f,g,k)~=~\tau^F(f,g,k)~=~\tau^F(h)~=~
\tau((h_0)^{-1}h)\\[1ex]
\hphantom{\tau^{iso}(f,g,k)~=~\tau^{iso}(h)}
\in W_1(A,\rho)~=~{\rm ker}(K_1(A_{\rho}[[z]])\to K_1(A))
\end{array}$$
with $F:A_{\rho}[[z]] \to A_{\rho}[[z]];z \mapsto 0$.

Continuing with our preparation for the construction of
$\widehat{B}^{-1}=\widehat{C}$, use the $\rho$-isomorphism
$$z^{-1}h_0~:~P \oplus z^kA^n \to P \oplus z^{-1}A^n~;~x \mapsto z^{-1}(h_0(x))$$
and the $A$-module isomorphisms
$$\begin{array}{l}
z^k(\theta_n)^k~:~A^n \to z^k A^n~;~ (a_1,a_2,\dots,a_n)
\mapsto z^k(\rho^k(a_1),\rho^k(a_2),\dots,\rho^k(a_n))\\[1ex]
z\theta_n~:~z^{-1}A^n \to A^n~;~
z^{-1}(a_1,a_2,\dots,a_n)\mapsto (\rho(a_1),\rho(a_2),\dots,\rho(a_n))
\end{array}$$
to define a $\rho$-isomorphism
$$\phi~=~(1 \oplus z\theta_n)(z^{-1}h_0)
(1 \oplus z^k(\theta_n)^k)~:~P \oplus A^n \to P \oplus A^n~.$$

We define an inverse for
$$\widehat C~=~\widehat C_1\oplus \widehat C_2\oplus \widehat C_3~:~
K_1(A,\rho) \oplus W_1(A,\rho) \oplus \widetilde{\mbox{\rm Nil}}_0(A,\rho^{-1}) \to K_1(A_{\rho}((z)))$$
by setting
$$\widehat B~=~\widehat B_1\oplus \widehat B_2\oplus \widehat B_3~:~
K_1(A_{\rho}((z))) \to K_1(A,\rho) \oplus W_1(A,\rho) \oplus \widetilde{\mbox{\rm Nil}}_0(A,\rho^{-1})$$
with
$$\begin{array}{l}
\widehat B_1~:~K_1(A_{\rho}((z)))\to K_1(A,\rho) ~;\\[1ex]
\hspace*{20mm} \tau(\alpha:A_{\rho}((z))^n \to A_{\rho}((z))^n) \mapsto
[P\oplus A^n,\phi]-[A^n,\theta_n]~,\\[1ex]
\widehat B_2~:~K_1(A_{\rho}((z)))\to W_1(A,\rho)  ~;\\[1ex]
\hspace*{20mm} \tau(\alpha:A_{\rho}((z))^n \to A_{\rho}((z))^n)
\mapsto w(f,g,k)~=~\tau((h_0)^{-1}h)~,\\[1ex]
\widehat B_3~:~K_1(A_{\rho}((z)))\to \widetilde{\mbox{\rm Nil}}_0(A,\rho^{-1})~;
\\[1ex]
\hspace*{20mm} \tau(\alpha:A_{\rho}((z))^n \to A_{\rho}((z))^n) \mapsto [P,\nu]~.
\end{array}
$$

We have to verify that the maps $\widehat B_i$ are well-defined, that
$$\widehat B_i \widehat C_j~=~\begin{cases} \Id &\hbox{if $i=j$} \\
0 &\hbox{if $i \neq j$}\end{cases}$$
and that
$$\widehat{C}_1\widehat{B}_1 +
\widehat{C}_2\widehat{B}_2 +\widehat{C}_3\widehat{B}_3~=~\Id~:~
K_1(A_{\rho}((z))) \to K_1(A_{\rho}((z))) ~.$$
The identities
$$\widehat B_i \widehat C_j~=~0~{\rm if}~i \neq j~~,~~
\widehat B_1\widehat C_1~=~\Id_{~K_1(A_{\rho}((z)))}~~,~~
\widehat B_3\widehat C_3~=~\Id_{ \widetilde{\mbox{\rm Nil}}_0(A,\rho^{-1})}$$
are easily reduced to the corresponding identites in the $K$-theory of Laurent polynomial ring
(subsection \ref{su:whlaurent}). The identity
$\widehat B_2\widehat C_2=\Id_{W_1(A,\rho)}$ is a matter of a trivial computation,
once we have proved that $\widehat B_2$ is well-defined. This
is proved below just after the next lemma, which proves
$\widehat C\widehat B=\Id$.

\bele\label{Lemma4} For every $A_{\rho}((z))$-module automorphism
$\alpha:A_{\rho}((z))^n \to A_{\rho}((z))^n$
$$\tau(\alpha)~=~(\widehat{C}_1\widehat{B}_1 +
\widehat{C}_2\widehat{B}_2 +\widehat{C}_3\widehat{B}_3)\tau(\alpha)
\in K_1(A_{\rho}((z)))~.$$
\enle
\noindent{\it Proof.} Apply $A_{\rho}((z))\otimes_{A_{\rho}[[z]]}-$ to the
chain equivalence of 1-dimensional $A$-induced f.g. projective $A_{\rho}[[z]]$-module
chain complexes
$$(f,g,k)~:~C(\widetilde\alpha:z^kA_{\rho}[[z]]^n \to A_{\rho}[[z]]^n)
\to C(z\zeta^{-1}-\nu:zP_\rho[[z]] \to P_{\rho}[[z]])$$
to obtain a chain equivalence of contractible
1-dimensional $A$-induced f.g. projective $A_{\rho}((z))$-module chain
complexes
$$1\otimes (f,g,k)~:~C(\widetilde\alpha:z^kA_{\rho}((z))^n \to A_{\rho}((z))^n)
\to C(z\zeta^{-1}-\nu:zP_\rho((z)) \to P_{\rho}((z)))~.$$
Use the $A_{\rho}((z))$-module isomorphisms
$$\aligned
&z^{-k}~:~z^kA_{\rho}((z))^n \to A_{\rho}((z))^n~,\\
&z^{-1}~:~zP_\rho((z)) \to P_{\rho}((z))
\endaligned$$
to define $A_{\rho}((z))$-module isomorphisms
$$\aligned
&C(\widetilde\alpha:z^kA_{\rho}((z))^n \to A_{\rho}((z))^n) ~\cong~
C(\alpha:A_{\rho}((z))^n \to A_{\rho}((z))^n) ~,\\
&C(z\zeta^{-1}-\nu:zP_\rho((z)) \to P_{\rho}((z)))~\cong~
C(1-z^{-1}\nu:P_\rho((z)) \to P_{\rho}((z)))
\endaligned$$
with
$$\widehat{C}_3\widehat{B}_3\tau(\alpha)~=~
\tau(1-z^{-1}\nu:P_\rho((z)) \to P_{\rho}((z))) \in K_1(A_{\rho}((z)))~.$$
It now follows from the sum formula for torsion that the morphism
$$\widehat{C}_2~:~W_1(A,\rho) \to K_1(A_{\rho}((z)))~;~w \mapsto w$$
sends $\widehat{B}_2\tau(\alpha)=w(f,g,k)$ to
$$\widehat{C}_2\widehat{B}_2\tau(\alpha)~
=~\tau(\alpha) - \widehat{C}_1\widehat{B}_1\tau(\alpha)-
\widehat{C}_3\widehat{B}_3\tau(\alpha)\in K_1(A_{\rho}((z)))~.\eqno{\qs}$$
\indent From now on, we shall write $\widehat B_2(\tau(\alpha))$ as
$\widehat B_2(\alpha)$. We have to check that $\widehat B_2$ is well defined,
that is
\begin{itemize}
\item[(i)] $\widehat B_2(\b\circ\a)=\widehat B_2(\b)+\widehat B_2(\a) \in W_1(A,\rho)$,\\[0.5ex]
\item[(ii)] $\widehat B_2(\alpha)=w(f,g,k)\in W_1(A,\rho)$ does not depend
on the choice of $k$.
\end{itemize}
\indent
The general observation is that the composite of two automorphisms of a f.g.
free $A_{\rho}((z))$-module
$$\alpha~=~\sum\limits^{\infty}_{j=-k}z^j\alpha_j~,~
\alpha'~=~\sum\limits^{\infty}_{j=-k'}z^j\alpha'_j~:~
A_{\rho}((z))^n \to A_{\rho}((z))^n$$
is an automorphism
$$\alpha''~=~\alpha'\alpha~=~\sum\limits^{\infty}_{j=-k''}z^j\alpha''_j~:~
A_{\rho}((z))^n \to A_{\rho}((z))^n~~(k''=k+k')$$
such that the cokernel $A_{\rho}[[z]]$-modules
$$\aligned
{}&P~=~{\rm coker}(\widetilde\alpha:z^kA_{\rho}[[z]]^n \to A_{\rho}[[z]]^n)~,\\
{}&P'~=~{\rm coker}(\widetilde\alpha':z^{k'}A_{\rho}[[z]]^n \to A_{\rho}[[z]]^n)~,\\
{}&P''~=~{\rm coker}(\widetilde\alpha'':z^{k''}A_{\rho}[[z]]^n \to A_{\rho}[[z]]^n)
\endaligned$$
fit into a short exact sequence
$$0 \to z^{k'}P \to P'' \to P' \to 0$$
which in principle gives (i). The case
$$\alpha'~=~1~:~A_{\rho}((z))^n \to A_{\rho}((z))^n~~,~~k'~=~1$$
in principle gives (ii). However, for the sake of precision we
shall now verify (i) and (ii) in detail.

We first make another general remark.
Let $P$ be an $A_{\rho}[[z]]$-module which admits an $A$-induced
f.g. projective $A_{\rho}[[z]]$-module resolution
$$R~:~0\rTo M_{\rho}[[z]]\rTo N_{\rho}[[z]] \rTo P \rTo 0~.$$
Any two such resolutions $R_1,R_2$ are related by a chain equivalence
$R_1 \to R_2$ with an isomorphism torsion
$\tau^{iso}(R_1 \to R_2) \in K^{iso}_1(A_{\rho}[[z]])$.
Denote the Witt vector torsion class (Definition \ref{Definition8}) by
$$\sigma(R_1,R_2)~=~w(\tau^{iso}(R_1 \to R_2)) \in W_1(A,\rho)~.$$
For any three such resolutions $R_1,R_2,R_3$ we have the sum formula
$$\sigma(R_1,R_3)~=~\sigma(R_1,R_2)+\sigma(R_2,R_3)~.$$
\indent Thus if $\a:\arzz^n\to\arzz^n$ is an automorphism
$$w(f,g,k)~=~\sigma(\mu(\a,k),\theta(\a,k))\in W_1(A,\rho)$$
where $\mu(\a,k)$ and $\theta(\a,k)$ are the following
resolutions of the module $P=P(k,\a)$:
\bq
\mu(\a,k)~:~0\rTo z^kA_{\rho}[[z]]^n\rTo^{\widetilde \alpha}
 A_{\rho}[[z]]^n \rTo^{\pi} P \rTo 0
\end{equation}
and
\bq
\theta(\a,k)~:~0 \rTo zP_\rho [[z]]   \rTo^{z\zeta^{-1}-\nu}
  P_{\rho}[[z]]   \rTo^{\widetilde \pi}  P \rTo  0 ~.
\end{equation}

\bepr\label{p:correc_k} The element
$$\sigma(\mu(\a,k),\theta(\a,k))\in W_1(A,\rho)$$
does not depend on $k$.\enpr
\noindent{\it Proof.} It suffices to check that the passage from $k$ to
$k+1$ does not change the invariant.  Let us consider the following
diagram

\begin{diagram}[LaTeXeqno]
 & & & 0  &   &  0 &  & 0 & & \\
& & & \dTo   &   &   \dTo &   & \dTo  \\
& & & z^{k+1}\arz^n  & \rTo^{\id}  & z^{k+1}\arz^n & \rInto & z^{k}\arz^n & & \\
& & & \dInto   & \boxed{\mu}   &   \dTo^{\a}  & \boxed{\mu\mu}  &   \dTo^{\a}
& &  \qquad\qquad \label{f:diag1}\\
& & & z^k \arz^n     &  \rInto   &   \arz^n        &  \rTo^{\id} &   \arz^n   & &   \\
& & & \dTo   &    &   \dTo &   &   \dTo & &  \\
&0& \rTo & P'    & \rTo^{i}    &   P(k+1,\a) & \rTo^p & P(k,\a) & \rTo & 0 & \\
\end{diagram}
\noindent
The two columns on the right are exactly $\mu(\a,k+1)$ and $\mu(\a,k)$.
The squares $\boxed{\mu}$ and $\boxed{\mu\mu}$ are obviously homotopy
commutative, and this implies the existence of the bottom line.  Recall
that the $A$-modules $P(k, \a), P(k+1,\a)$ come equipped with nilpotent
$\rho^{-1}$-homomorphisms $\nu_k$ and, respectively, $\nu_{k+1}$.
Endow the $A$-module $P'\cong z^k A$ with the zero nilpotent
endomorphism, then the bottom line becomes an exact sequence in the
category of modules and their nilpotent $\rho^{-1}$-endomorphisms.
This exact sequence induces the following diagram:

\begin{diagram}[LaTeXeqno]
 & & & 0  &   &  0 &  & 0 & & \\
& & & \dTo   &   &   \dTo &   & \dTo  \\
& & & zP'_\rho[[z]]  & \rTo^{\id}  &  zP(k+1,\a)_\rho[[z]] & \rTo & zP(k,\a)_\rho[[z]] & & \\
& & & \dInto   & \boxed{\vartheta}   &   \dTo  & \boxed{\vartheta\vartheta}  &\dTo & &
\qquad\qquad\label{f:diag2}
\\
& & & P'_\rho[[z]]  & \rTo^{\id}  &  P(k+1,\a)_\rho[[z]] & \rTo & P(k,\a)_\rho[[z]] & & \\
& & & \dTo   &    &   \dTo &   &   \dTo & &  \\
&0&\rTo & P'    & \rTo^{i}    &   P(k+1,\a) & \rTo^p & P(k,\a) & \rTo & 0 & \\
\end{diagram}
\noindent where all the rows are exact sequences, and the two columns
on the right are the resolutions $\theta(k+1,\a)$ and resp.
$\theta(k,\a)$.  The computation of $\sigma(\mu(\a,k+1),\theta(\a,k+1))$ is
reduced to the comparison of the resolutions represented by the central
columns of the both diagrams.

Consider the square $\boxed{\mu}$
as a map of two resolutions corresponding to
modules $P'$ and $P(k+1,\a)$.
Using  the diagram
(\ref{f:diag1})
it is easy to prove:
\bele\label{l:kappa}
There is an epimorphic chain homotopy equivalence
$$C(\boxed{\mu})\to \mu(\a,k)$$
with kernel the chain complex
\bq
\{0\rTo z^{k+1}\arz^n\rTo^{\id} z^{k+1}\arz^n\rTo 0\}
\end{equation}
\enle
\noindent {\it Sketch of the proof.}
The composition of two maps of resolutions, - the first corresponding
to $\boxed{\mu}$, the second to $\boxed{\mu\mu}$ is homotopic to zero,
being a lifting of the zero map $P'\to P(k,\a)$.  This implies a direct
construction of the map $C(\boxed{\mu})\to \mu(\a,k)$, and checking
through this construction leads to the proof of the lemma.  $\qs$

Let
\bq
\xi= \{0\rTo zP'[[z]] \rTo  P'[[z]]\rTo 0 \}
\end{equation}
Similar reasoning gives the following:

\bele
\label{l:Q}
There is an epimorphic $A$-induced chain homotopy equivalence
$\pi:C(\boxed{\vartheta})\to \theta(\a,k)$
with kernel
$${\rm ker}(\pi)~=~C(\id:\xi\to\xi)~.\eqno{\qs}$$
\enle

In particular $C(\boxed{\vartheta})$ is a resolution for $S(k,\a)$.
It follows from the preceding lemmas that
$$\sigma(C(\boxed{\vartheta}),\theta(\a,k))~=~0~=~\sigma(C(\boxed{\mu}),
\mu(\a,k))~,$$
and therefore
\bq
\sigma(\mu(\a,k),\theta(\a,k))~=~\sigma(C(\boxed{\mu}), C(\boxed{\vartheta}))
\end{equation}

Let
\bq
\eta= \{0\rTo z^{k+1} \arz^n \rTo z^k\arz^n \rTo 0 \}
\end{equation}
\bepr
There is a commutative diagram of $\arz$-module chain complexes
and chain maps
\begin{diagram}[LaTeXeqno]
0& \rTo  & \mu(\a,k+1) & \rTo & C(\boxed{\mu}) & \rTo & \Sigma(\eta) & \rTo & 0 \\
 &       & \dTo        &      & \dTo            &       & \dTo         &      &   \quad\label{f:map_cone} \\
0& \rTo  & \theta(\a,k+1) & \rTo & C(\boxed{\vartheta}) & \rTo & \Sigma(\xi) & \rTo & 0 \\
\end{diagram}
\noindent (here $\Sigma$ stands for suspension in the category of chain
complexes, and all the horizontal arrows are $A$-induced
homomorphisms).\enpr
\noindent{\it Proof.} Consider the following square of chain complexes.
\begin{diagram}[size=2em,LaTeXeqno]
\eta & \rTo & \mu(k+1,\a)\\
\dTo &   &\dTo   \qquad\label{f:squsqu}\\
\xi    & \rTo & \theta(k+1,\a)\\
\end{diagram}

\noindent
(Here the upper horizontal arrow comes from the square $\boxed{\mu}$,
the lower comes from $\boxed{\vartheta}$, the vertical arrows are the
maps of resolutions, induced by the identity maps of the corresponding
complexes.) The square (\ref{f:squsqu}) is homotopy commutative
(indeed, all these complexes are actually resolutions of the
corresponding $\arz$-modules, and to check the homotopy commutativity it
suffices to check the commutativity on the level of the modules
themselves, which is obvious).  By the basic properties of the
algebraic mapping cone construction this implies the existence of the
left square (strictly commutative) of (\ref{f:map_cone}) and the rest
of this diagram is now obvious.  $\qs$

Therefore
\bq
\sigma(C(\boxed{\mu}), C(\boxed{\theta}))~=~
\sigma(\mu(\a,k+1),\theta(\a,k+1))+\sigma(\Sigma(\eta),\Sigma(\xi))
\end{equation}
and the proof of Proposition \ref{p:correc_k} is concluded.\hfill$\qs$

We can now identify
$$\widehat B_2(\a)~=~\sigma(\mu(\a,k),\theta(\a,k)) \in W_1(A,\rho)$$
and  proceed to the next step of the verification.
It is quite obvious that
$$\widehat B_2(\a\oplus\b)~=~\widehat B_2(\a)+\widehat B_2(\b)~~,~~
\widehat B_2(\id)~=~0 \in W_1(A,\rho)~,$$
so $\widehat B_2$ is a well defined map $\karzz\to
W_1(A,\rho)$. It remains to check that this map is a homomorphism.

\bepr\label{p:group_hom}
$\widehat B_2(\b\circ\a)= \widehat B_2(\b)+\widehat B_2(\a) \in W_1(A,\rho).$
\enpr

\noindent{\it Proof.}
Note first of all that it suffices to prove the proposition for the
particular case when $\b(\arz^n)\subset \arz^n$.  (Indeed, the general case
is easily reduced to this particular one considering $\b'=\zeta^N\b$
where $\zeta: \arz^n\to \arz^n$ is defined by
$\zeta(x_1,\dots,x_n)=(z\rho(x_1),\dots,z\rho(x_n))$, and $N$ is
sufficiently large.)

For this particular case the proof goes on the same lines as the proof
for Proposition \ref{p:correc_k}.  Consider the following diagram
\begin{diagram}[LaTeXeqno]
 & & & 0  &   &  0 &  & 0 & & \\
& & & \dTo   &   &   \dTo &   & \dTo  \\
& & & z^k \arz^n  & \rTo^{\id}  & z^k\arz^n  & \rTo^\a & z^k\arz^n & & \\
& & & \dTo^\a   &   &   \dTo^{\b\a}  &   &   \dTo^{\b} & &  \qquad\qquad \label{f:diaggg}\\
& & & \arz^n    &  \rTo^\b   &   \arz^n        &  \rTo^{\id} &   \arz^n   & &   \\
& & & \dTo   &    &   \dTo &   &   \dTo & &  \\
&0& \rTo & P(k,\a)    & \rTo^{i}    &   P(k,\b\a) & \rTo^p & P(k,\b) & \rTo & 0 & \\
\end{diagram}
\noindent where the squares  in the middle are obvious, and the
bottom line is obtained from these two squares.
Form the diagram
\begin{diagram}[LaTeXeqno]
 & & & 0  &   &  0 &  & 0 & & \\
& & & \dTo   &   &   \dTo &   & \dTo  \\
&0 &\rTo & zP(k,\a)[[t]]  & \rTo  &  zP(k,\b\a)[[t]] & \rTo & zP(k,\b)[[t]] &\rTo &0 \\
& & & \dTo   &    &   \dTo  &   &\dTo & &
\qquad\qquad\label{f:diaggggggg}
\\
&0 &\rTo & P(k,\a)[[z]]  & \rTo  &  P(k,\b\a)[[z]] & \rTo & P(k,\b)[[z]] &\rTo &0 \\
& & & \dTo   &    &   \dTo &   &   \dTo & &  \\
&0&\rTo & P(k,\a)    & \rTo^{i}    &   P(k,\b\a) & \rTo^p & P(k,\b) & \rTo & 0 & \\
\end{diagram}

The rest of the proof repeats step by step the proof of Proposition \ref{p:correc_k}
with the corresponding changes.
$\qs$

\section{On the image of the Witt vectors in the Whitehead group}
\label{s:ex}

The Main Theorem reduces the computation of the group
$K_1(A_\rho((z)))$ to the computation of three groups.  Two of these
three summands are classical algebraic $K$-theoretic groups\,:
$K_1(A,\rho)$ and $\widetilde{\Nil}_0(A,\rho^{-1})$, which have already
appeared in the computation of $K_1(A_{\rho}[z,z^{-1}])$.  Much less is
known about the third summand $\woar$.  By definition
$$\woar~=~{\rm im}\bigg(\war\rTo K_1(A_{\rho}((z)))\bigg)$$
with
$$\war~=~\{1+a_1z+a_2z^2+\dots \,\vert\,a_i \in A\} \subseteq \gA$$
the subgroup of the Witt vectors in the multiplicative
group of units $\gA=A_{\rho}[[z]]^{\bullet}$.
The aim of this section is to obtain some partial information
about the abelian group $\woar$.

For any group $G$ we write the commutators in the usual fashion as
$$[g,g']~=~gg'g^{-1}g^{\prime -1} \in G~~(g,g' \in G)~,$$
and denote by $G^{ab}$ the abelian quotient of $G$ by the
normal subgroup $[G,G] \subseteq G$ of the elements $x \in G$
which can be expressed as products of commutators
$$x~=~\prod\limits^n_{r=1}[g_r,g'_r] \in G~.$$
There is a natural surjection
\bq J~:~\warab \to \woar
\end{equation}
so that $\warab$ is in a sense a first approximation to $\woar$, and
the problem of computing of $\woar$ may be viewed as the problem of
computing the group ${\rm ker}\, J$.  The next step of approximation can be
obtained via $\gA$ and the factorization
\begin{diagram}
\warab & \rTo^{J'} & \gA^{ab} \\
       &  \rdTo(2,2)^{J}& \dTo   \\
       &              & K_1(A_{\rho}((z))) \\
\end{diagram}
where the vertical arrow is induced by the natural inclusion.  In
particular, ${\rm ker}\, J'\subset {\rm ker}\, J$.  It turns out that ${\rm ker}\, J'$ is
easily described in terms of the group $\warab$ and the conjugation
action of $\gA$ on this group.  This can be done in a more general
framework of semidirect product of groups, and this is the aim of the
next section.

\subsection{On semidirect products of groups}
\label{su:absgroup}

The semidirect product of groups $H,K$ twisted by a homomorphism
$\xi:K \to \text{Aut}(H)$ is the group
$$G~=~H\times_{\xi}K$$
in which every element $x \in G$ has a unique expression as a product
$$x~=~hk \in G~~(h \in H,k \in K)$$
with
$$(h_1k_1)(h_2k_2)~=~h_1\xi(k_1)(h_2)k_1k_2 \in G~~,~~
\xi(k)(h)~=~khk^{-1} \in H$$
and there is defined a short exact sequence of groups
\bq
1\rTo H \rTo^{i} G \rTo^{j} K \rTo 1\end{equation}
with
$$i~:~H \to G~;~h \mapsto h~,~j~:~G \to K~;~hk \mapsto k~.$$
The abelianization $H^{ab}$ is a $\ZZZ[K]$-module via
$$\ZZZ[K] \times H^{ab} \to H^{ab}~;~
(nk,h) \mapsto \xi(k)(h^n)~=~kh^nk^{-1}~.$$

In our applications we shall have
$$G~=~A_{\rho}[[z]]^{\bullet}~=~\gA~=~H\times_{\xi}K$$
with
$$H~=~\war~,~K~=~A^{\bullet}~,~\xi:K \to \text{Aut}(H)~;~
k \mapsto (h \mapsto khk^{-1})~.$$

\bepr\label{p:kergr} For any semidirect product $G=H\times_{\xi}K$
the kernel of the induced morphism of
abelian groups $i^{ab}:H^{ab} \to G^{ab}$ is given by
$${\rm ker}(i^{ab})~=~{\rm ker}(\epsilon)H^{ab} \subseteq H^{ab}$$
with $\epsilon:\ZZZ[K]\to \ZZZ;k \mapsto 0$ the usual augmentation.\enpr
\noindent{\it Proof.} (i) We prove first that ${\rm ker}\, i^{ab}
\supset ({\rm ker}\,\epsilon)(H^{ab})$.
It suffices to observe that for any $h\in H$, $k\in K$
$$i(\xi(k)(h)-h)~=~[k,h] \in G~,$$
so that
$$i^{ab}(\xi(k)(h)-h)~=~0 \in G^{ab}~.$$
(ii) Now for the reverse inclusion\,: we have to show that if $x \in H$
is such that $i(x) \in [G,G]$ then $[x]=y(h)\in H^{ab}$ for some
$y \in {\rm ker}(\epsilon)$, $h \in H$. For any $hk,h'k' \in G$ we have
$$[hk,h'k']~=~\xi(k)(h)\xi(k'k^{-1})(hh')\xi(k^{\prime -1})(hh'h^{-1})
[h,h'][k,k'] \in [G,G]~.$$
If $x \in H$ is such that $i(x) \in [G,G]$ then
$$x~=~\prod^n_{r=1}[h_rk_r,h'_rk'_r]  \in G$$
for some $h_r \in H$, $k_r \in K$ with
$$\prod^n_{r=1}[k_r,k'_r]~=~1 \in K$$
and
$$\begin{array}{l}
i(x)~=~\sum\limits^n_{r=1}
\xi(k)(h)\xi(k_r'k_r^{-1})(h_rh_r')\xi(k_r^{\prime -1})
(h_rh_r'h_r^{-1})\\[1ex]
\hspace*{15mm}
=~\sum\limits^n_{r=1}\bigg((\xi(k_r)-\xi(k_r^{\prime -1}))(h_r) +
(\xi(k_r'k_r^{-1})-\xi(k_r^{\prime -1}))(h_rh_r')\bigg)\\[1ex]
\hspace*{30mm}\in {\rm ker}\, (\epsilon)H^{ab} \subseteq H^{ab}~.
\end{array}$$
$\qs$

\subsection{Explicit examples }
\label{su:expex}

Now we return to Witt vectors. The direct application of
Proposition \ref{p:kergr} gives
\bq
{\rm ker}\, J' = ({\rm ker}\,\epsilon)(\warab)
\end{equation}
where $\epsilon:\ZZZ\gA\to\ZZZ$ is the augmentation in the group ring of the
group $\gA$.  This result leads quickly to the construction of
non-trivial elements in ${\rm ker}\, J'$.  Indeed, any element of type
$y=\xi(\alpha)(x)-x$ where $\a\in\gA, x\in\war$ is in ${\rm ker}\, J'$.  If
$$x~=~1+\b z~~(\beta\in A)$$
then $y$ has a following representative:
\bq
\bar y~=~\a(1+\b z)\a^{-1} (1+\b z)^{-1}
\in
\war
\end{equation}
Assume now that $\rho=\id$, and $\a\beta\a^{-1}\not=\b$.  The element
$\bar y$ is of the form \bq 1+(\a\beta\a^{-1} - \beta)z +O(z^2)
\end{equation} and it does not belong to $[\war, \war]$ since every
element
$$1+a_1z+a_2z^2+\dots \in [\war, \war]$$
obviously has $a_1=0$.

\label{refer}

\end{document}